\documentclass{amsart}
\usepackage{amssymb}

\usepackage{amscd}

\newtheorem{theorem}{Theorem}[section]

\newtheorem{lemma}[theorem]{Lemma}
\newtheorem{proposition}[theorem]{Proposition}

\theoremstyle{definition}
\newtheorem{definition}{Definition}[section]

\theoremstyle{remark}

\begin{document}
\title{A local proof of Petri's conjecture at the general curve}
\author{Herb Clemens}
\date{May, 2006}
\curraddr{Mathematics Department, Ohio State University, Columbus, OH 43210, USA}
\email{clemens@math.ohio-state.edu}
\maketitle

\begin{abstract}
A proof of Petri's general conjecture on the unobstructedness of linear
systems on a general curve is given, using only the local properties of the
deformation space of the pair (curve, line bundle).
\end{abstract}

\section{Introduction\label{0}}

Let $L_{0}$ denote a holomorphic line bundle of degree $d$ over a compact
Riemann surface $C_{0}$. The Petri conjecture stated that, if $C_{0}$ is a
curve of general moduli, the mapping 
\begin{equation*}
\mu _{0}:H^{0}\left( L_{0}\right) \otimes H^{0}\left( \omega _{C_{0}}\otimes
L_{0}^{\vee }\right) \rightarrow H^{0}\left( \omega _{C_{0}}\right) ,
\end{equation*}
is injective. Later, this assertion was given a more modern interpretation
making it a central question in the study of curves and their linear
series---what is now called Brill-Noether theory.

To recap the modern formulation we proceed as in \cite{AC}. Let $%
C_{0}^{\left( d\right) }$ denote the $d$-th symmetric product of $C_{0}$ and
let $\Delta \subseteq C_{0}^{\left( d\right) }\times C_{0}$ denote the
tautological divisor. Let 
\begin{equation*}
\Bbb{P}^{r}=\Bbb{P}\left( H^{0}\left( L_{0}\right) \right) .
\end{equation*}
For the projection 
\begin{equation*}
p_{*}:C_{0}^{\left( d\right) }\times C_{0}\rightarrow C_{0}^{\left( d\right)
}
\end{equation*}
and exact sequence 
\begin{equation*}
0\rightarrow O_{C_{0}^{\left( d\right) }\times C_{0}}\rightarrow
O_{C_{0}^{\left( d\right) }\times C_{0}}\left( \Delta \right) \rightarrow
\left. O_{C_{0}^{\left( d\right) }\times C_{0}}\left( \Delta \right) \right|
_{\Delta }\rightarrow 0,
\end{equation*}
one has that 
\begin{equation*}
T_{C_{0}^{\left( d\right) }}=p_{*}\left( \left. O_{C_{0}^{\left( d\right)
}\times C_{0}}\left( \Delta \right) \right| _{\Delta }\right) .
\end{equation*}
Applying the derived functor $Rp_{*}\circ \mathcal{O}_{\Bbb{P}^{r}\times
C_{0}}$ to the above exact sequence as in (2.6) of \cite{AC}, one obtains an
exact sequence 
\begin{equation*}
0\rightarrow N_{\Bbb{P}^{r}\backslash C_{0}^{\left( d\right) }}\rightarrow
O_{\Bbb{P}^{r}}\otimes H^{1}\left( O_{C_{0}}\right) \rightarrow O_{\Bbb{P}%
^{r}}\left( 1\right) \otimes H^{1}\left( L_{0}\right) \rightarrow 0,
\end{equation*}
where $N_{A\backslash B}$ denotes the normal bundle of $A$ in $B$. So the
dual of the kernel of $\mu _{0}$ above is exactly 
\begin{equation*}
H^{1}\left( \left. N_{\Bbb{P}\backslash C_{0}^{\left( d\right) }}\right| _{%
\Bbb{P}^{r}}\right) .
\end{equation*}

Via the standard short exact sequence of normal bundles, Petri's conjecture
becomes the assertion 
\begin{equation*}
H^{1}\left( N_{\Bbb{P}^{r}\backslash C_{0}^{\left( d\right) }}\right) =0,
\end{equation*}
that is, the deformation theory of linear series is unobstructed at a curve
of general moduli.

There are several proofs of Petri's conjecture, proofs via degeneration by
Gieseker \cite{Gi} and Eisenbud-Harris \cite{EH} and a proof via
specialization to the locus of curves on a general $K3$-surface due to
Lazarsfeld \cite{L} (see also \cite{P}). However the only proof based on
properties of the infinitesimal deformation of the general curve, as opposed
to some specialization of it, is a proof for $r\leq 2$ by Arbarello and
Cornalba in \cite{AC}. In conversations concerning his joint work with
Cornalba, Arbarello explained to the author the viewpoint of \cite{ACGH}
that there should exist a generalization to higher order of the following
result (which appears both in \cite{ACGH} and \cite{AC}):

Let 
\begin{equation*}
\frak{D}_{n}\left( L_{0}\right)
\end{equation*}
denote the sheaf of holomorphic differential operators of order $\leq n$ on
sections of the line bundle $L_{0}$. (If 
\begin{equation*}
L_{0}=\mathcal{O}_{C_{0}}
\end{equation*}
we denote this sheaf simply as $\frak{D}_{n}$.) The first-order deformations
the pair $\left( L_{0},C_{0}\right) $ are in natural one-to-one
correspondence with the elements 
\begin{equation*}
\psi \in H^{1}\left( \frak{D}_{1}\left( L_{0}\right) \right)
\end{equation*}
in such a way that a section $s_{0}$ of $L_{0}$ deforms to first order with
the deformation $\psi $ if and only if the element 
\begin{equation*}
\psi \left( s_{0}\right) \in H^{1}\left( L_{0}\right)
\end{equation*}
is zero.

Furthermore he pointed out that an appropriate higher-order generalization
of this fact and a simple Wronskian argument would immediately yield a
``local'' proof of Petri's general conjecture at the general curve (see \S 
\ref{3} below). The purpose of this paper is to carry out that
generalization.

The general idea of the proof is to use the Kuranishi theory of
(curvilinear) $C^{\infty }$-trivializations of deformations of complex
manifolds as it applies to the total space the dual line bundle $L_{0}^{\vee
}$. Roughly speaking, if we denote the $t$-disk as $\Delta $ and are given a 
$C^{\infty }$-trivialization 
\begin{equation*}
F_{\sigma }=\left( \sigma ,\pi \right) :M\rightarrow M_{0}\times \Delta
\end{equation*}
of a deformation $M/\Delta $ of a complex manifold $M_{0}$, Kuranishi
associated to this situation a power series 
\begin{equation*}
\xi =\xi _{1}t+\xi _{2}t^{2}+\ldots
\end{equation*}
where each $\xi _{j}$ is a $\left( 0,1\right) $-form with coefficients in (a
subsheaf of) the tangent bundle of $M_{0}$. $F_{\sigma }$ is not allowed to
be an arbitrary $C^{\infty }$-isomorphism over $\Delta $. The relevant
restriction is that trajectory of each point on $M_{0}$ must be holomorphic,
that is, 
\begin{equation*}
\sigma ^{-1}\left( x_{0}\right) \subseteq M
\end{equation*}
must be a holomorphic disk for each $x_{0}\in M_{0}$. This is of course just
a restriction on the choice of trivialization; it implies no restriction on
the deformation $M/\Delta $. For such a trivialization, the holomorphic
functions $f$ on $M$ have a very nice form; namely we can write power-series
expansions 
\begin{equation*}
f\circ F_{\sigma }^{-1}=f_{0}+f_{1}t+f_{2}t^{2}+\ldots
\end{equation*}
such that the holomorphicity condition 
\begin{equation*}
\overline{\partial }_{M}f=0
\end{equation*}
becomes just 
\begin{equation*}
\left( \overline{\partial }_{M_{0}}-\xi \right) \left(
f_{0}+f_{1}t+f_{2}t^{2}+\ldots \right) =0.
\end{equation*}

Although later on we will actually need to consider a slightly more general
case in the body of this paper, it is perhaps helpful as an introduction to
give the line of reasoning of the paper in the case in which $M_{0}$ happens
to be the total space of a holomorphic line bundle 
\begin{equation*}
q_{0}:L_{0}^{\vee }\rightarrow C_{0}
\end{equation*}
over a compact Riemann surface $C_{0}.$ One easily sees that the deformation
is a deformation of holomorphic line bundles if and only if the Kuranishi
data $\xi ^{L}$ are invariant under the action of the $\Bbb{C}^{*}$-action
on $L_{0}^{\vee }$. In fact, if $\chi $ denotes the $\left( 1,0\right) $
Euler vector field on $L_{0}^{\vee }$ associated with the natural $\Bbb{C}
^{*}$-action on the line bundle, this is just the condition 
\begin{equation*}
\left[ \chi ,\xi _{j}^{L}\right] =0
\end{equation*}
for all $j$, that is, that the $\xi _{j}^{L}$ can be written everywhere
locally in the form 
\begin{equation}
q_{0}^{*}\left( \alpha \right) \cdot \chi +q_{0}^{*}\left( \beta \right)
\cdot \tau _{L}  \label{A}
\end{equation}
where $\alpha $ and $\beta $ are $\left( 0,1\right) $-forms on $C_{0}$ and $%
\tau _{L}$ is a lifting of a $\left( 1,0\right) $-vector-field $\tau _{C}$
on $C_{0}$ such that 
\begin{equation*}
\left[ \chi ,\tau _{L}\right] =0.
\end{equation*}
(The ``associated'' or ``compatible'' Kuranishi data for the deformation of $%
C_{0}$ is just given by $\xi _{j}^{C}=\beta \cdot \tau _{C}$.) Sections $s$
of $L$ are just functions $f$ on $L_{0}^{\vee }$ for which 
\begin{equation*}
L_{\chi }\left( f\right) =f
\end{equation*}
where $L_{\chi }$ denotes Lie differentiation with respect to the vector
field $\chi $.

Suppose now we have a line-bundle deformation $\left( L/\Delta ,C/\Delta
\right) $ of $\left( L_{0},C_{0}\right) $ with compatible trivializations 
\begin{eqnarray*}
\sigma &:&C\rightarrow C_{0} \\
\lambda &:&L^{\vee }\rightarrow L_{0}^{\vee }
\end{eqnarray*}
and a section $s$ of $L$ whose zeros are given by 
\begin{equation*}
\sigma ^{-1}\left( zeros\left( s_{0}\right) \right) .
\end{equation*}
Rescaling $\lambda $ in the fiber direction we arrive at a trivialization of
the deformation $L^{\vee }$ of $L_{0}^{\vee }$ for which $s$ is constant,
that is, 
\begin{equation*}
s=s_{0}\circ \lambda .
\end{equation*}
We call such compatible trivializations of $C_{0}$ and $L_{0}^{\vee }$
``adapted'' to the section $s$.

Of course we have twisted the almost complex structure on $C_{0}$ and $%
L_{0}^{\vee }$ to achieve this trivialization. To keep track of this
twisting, we consider only ``Schiffer-type'' deformations $C$ of $C_{0}$,
for which the twist in almost complex structure is given almost everywhere
by a gauge transformation, that is, by a power series 
\begin{equation*}
\beta =\beta _{1}t+\beta _{2}t^{2}+\ldots
\end{equation*}
where the $\beta _{j}$ are $C^{\infty }$-vector-fields of type $\left(
1,0\right) $ on $C_{0}-\left\{ p\right\} $ and meromorphic in a small
analytic neighborhood of $p$. Then we take 
\begin{equation*}
\xi ^{C}=\frac{e^{\left[ \beta ,\ \right] }-1}{\left[ \beta ,\ \right] }%
\left( \overline{\partial }_{C_{0}}\beta \right)
\end{equation*}
(see \cite{GM}) and get a compatible trivialization of $L^{\vee }/\Delta $
by lifting the $\beta _{j}$ to vector fields $\tilde{\beta}_{j}$ on $%
L_{0}^{\vee }$ with 
\begin{equation*}
\left[ \tilde{\beta}_{j},\chi \right] =0
\end{equation*}
with the same meromorphic property near $q_{0}^{-1}\left( p\right) $.
Holomorphicity of a section $s$ becomes the condition 
\begin{equation*}
\left( \overline{\partial }_{L_{0}^{\vee }}\left( e^{L_{-\beta }}\left(
f\right) \right) \right) =0
\end{equation*}
on the power series 
\begin{equation*}
f=f_{0}+f_{1}t+f_{2}t^{2}+\ldots
\end{equation*}
representing $s$ as a function on $L_{0}^{\vee }\times \Delta $. That is,
the condition is simply that the pull-back of $f$ via the gauge
transformation is a power series whose coefficients are meromorphic sections
of $L_{0}$.

If we have a holomorphic section $s$ of $L$ whose restriction to $s_{0}$ has
simple zeros $D_{0}$ and if $\tilde{\beta}$ is zero in a small analytic
neighborhood of 
\begin{equation*}
D=zero\left( s\right) \subset C
\end{equation*}
then there is a $C^{\infty }$-automorphism 
\begin{equation*}
\Phi :C_{0}\times \Delta \rightarrow C_{0}\times \Delta
\end{equation*}
defined over $\Delta $ such that:

\begin{enumerate}
\item  $\Phi $ is holomorphic in a small analytic neighborhood of $D\cup
\left\{ p\right\} $.

\item  
\begin{equation*}
\Phi \left( \left\{ x_{0}\right\} \times \Delta \right) 
\end{equation*}
is a holomorphic disk for each $x_{0}\in C_{0}$.

\item  
\begin{equation*}
\Phi \circ F_{\sigma }\left( D\right) =D_{0}\times \Delta .
\end{equation*}
\end{enumerate}

The rough (imprecise) idea is that trivialization $\Phi \circ F_{\sigma }$
can also be considered to be of Schiffer type for some vector field 
\begin{equation*}
\gamma =\gamma _{1}t+\gamma _{2}t^{2}+\ldots .
\end{equation*}
$\gamma $ lifts to a vector field $\tilde{\gamma}$ associated to a
Schiffer-type trivialization of the deformation $L^{\vee }/\Delta $ of $%
L_{0}^{\vee }$ which is adapted to the section $s$. Since by construction $s$
corresponds to the ``constant'' power series 
\begin{equation*}
f_{0}+0\cdot t+0\cdot t^{2}+\ldots ,
\end{equation*}
we have the equation 
\begin{equation*}
\left( \overline{\partial }_{L_{0}^{\vee }}\left( e^{L_{-\tilde{\gamma}%
}}\left( f_{0}\right) \right) \right) =0
\end{equation*}
that is, 
\begin{equation}
\left[ \overline{\partial }_{L_{0}^{\vee }},e^{L_{-\tilde{\gamma}}}\right]
\left( f_{0}\right) =0.  \label{B}
\end{equation}
It is in this way that we produce elements of $H^{1}\left( \frak{D}%
_{n+1}\left( L_{0}\right) \right) $ for all $n\geq 0$ which must annihilate
sections $s_{0}$ of $L_{0}$ which extend to sections of $L$. (The difficulty
is of course that the elements of $H^{1}\left( \frak{D}_{n+1}\left(
L_{0}\right) \right) $ depend on the choice of $s_{0}$. To remedy this we
will eventually have to replace the deformation $C/\Delta $ of $C_{0}$ with
the deformation 
\begin{equation*}
\Bbb{P}/\Delta =\Bbb{P}\left( H^{0}\left( L/\Delta \right) \right)
\end{equation*}
of $\Bbb{P}\left( H^{0}\left( L_{0}\right) \right) $ and replace $L$ with $%
\mathcal{O}\left( 1\right) $.)

As one of the simplest concrete examples, let 
\begin{equation*}
C_{0}=\frac{\Bbb{C}}{\Bbb{Z}+\Bbb{Z}\sqrt{-1}}
\end{equation*}
with linear holomorphic coordinate $z$ on $\Bbb{C}$. For a $C^{\infty }$%
-function $\rho $ supported on $\left\{ z:\left| z\right| \leq 1/8\right\} $
and identically $1$ on $\left\{ z:\left| z\right| \leq 1/16\right\} $, let 
\begin{eqnarray*}
\beta _{1} &=&\frac{\rho }{z}\cdot \frac{\partial }{\partial z} \\
\beta _{j} &=&0,\ j>1.
\end{eqnarray*}
This is a non-trivial deformation since, to first order it is given by the
generator 
\begin{equation*}
\overline{\partial }_{C_{0}}\left( \frac{\rho }{z}\cdot \frac{\partial }{%
\partial z}\right) \in H^{1}\left( T_{C_{0}}\right) .
\end{equation*}

For $L_{0}$ we can take the line bundle of degree $2$ given by the divisor 
\begin{equation*}
D_{0}=\left\{ \frac{1+\sqrt{-1}}{4}\right\} +\left\{ \frac{3+3\sqrt{-1}}{4}
\right\} .
\end{equation*}
with corresponding section $s_{0}$. Let $s$ be some extension of the section 
$s_{0}$. For a trivialization 
\begin{equation*}
F_{\sigma }:C\rightarrow C_{0}\times \Delta
\end{equation*}
associated to the above Kuranishi data, the zero set $D=D^{\prime
}+D^{\prime \prime }$ of the section $s$ is given by two power series 
\begin{eqnarray*}
z &=&a\left( t\right) =\frac{1+\sqrt{-1}}{4}+a_{1}t+\ldots \\
z &=&b\left( t\right) =\frac{3+3\sqrt{-1}}{4}+b_{1}t+\ldots
\end{eqnarray*}
since the deformation of (almost) complex structure is zero near $D_{0}$. So
near $D^{\prime }$ we recursively solve for 
\begin{equation*}
\Phi \left( z,t\right) =\left( a^{\prime }\left( z,t\right) ,t\right)
\end{equation*}
such that 
\begin{equation*}
a^{\prime }\left( a\left( t\right) ,t\right) \equiv \frac{1+\sqrt{-1}}{4}
\end{equation*}
and similarly near $D^{\prime \prime }$ for 
\begin{equation*}
\Phi \left( z,t\right) =\left( b^{\prime }\left( z,t\right) ,t\right)
\end{equation*}
such that 
\begin{equation*}
b^{\prime }\left( b\left( t\right) ,t\right) \equiv \frac{3+3\sqrt{-1}}{4}.
\end{equation*}
Near $\left\{ 0\right\} \times \Delta $ we take 
\begin{equation*}
\Phi \left( z,t\right) =\left( z,t\right)
\end{equation*}
and then extend $\Phi $ to a family of diffeomorphism on all of $C_{0}$ by a 
$C^{\infty }$ patching argument. For the new trivialization 
\begin{equation*}
\Phi \circ F_{\sigma }:C\rightarrow C_{0}\times \Delta
\end{equation*}
the divisor $D$ giving the line bundle $L$ is ``constant'' so that the
pull-back of $s_{0}$ via the product structure gives rise to a compatible
trivialization of $L$.

The Petri proof will follow from doing this process (for a line-bundle
deformation of $L_{0}$ for which all sections extend) for every
Schiffer-type variation of a generic curve $C_{0}$. We show that the set of
equations $\left( \ref{B}\right) $ we obtain implies that the higher $\mu $%
-maps 
\begin{equation*}
\mu _{n+1}:\ker \left( \mu _{n}\right) \rightarrow H^{0}\left( \omega
_{C_{0}}^{n+2}\right) =H^{1}\left( T_{C_{0}}^{n+1}\right)
\end{equation*}
are all zero. As Arbarello-Cornalba-Griffiths-Harris showed twenty years
ago, this implies Petri's conjecture.

We shall use Dolbeault cohomology throughout this paper. In particular, the
sheaf $\frak{D}_{n}\left( L_{0}\right) $ has both both a left and a right $%
\mathcal{O}_{C_{0}}$-module structure and we define 
\begin{equation*}
A^{0,i}\left( \frak{D}_{n}\left( L_{0}\right) \right)
:=A_{C_{0}}^{0,i}\otimes _{\mathcal{O}_{C_{0}}}\frak{D}_{n}\left(
L_{0}\right)
\end{equation*}
where $A^{0,i}$ is the sheaf of $C^{\infty }$-$\left( 0,i\right) $-forms.
Also the context will hopefully eliminate any confusion between two standard
notations used in this paper, namely the notation $L$ and $L_{0}$ for line
bundles and the notation 
\begin{equation*}
L_{\tau }^{k}=\underset{k-times}{\underbrace{L_{\tau }\circ \ldots \circ
L_{\tau }}}
\end{equation*}
where $L_{\tau }$ denotes Lie differentiation with respect to a vector field 
$\tau $.

The author wishes to thank E. Arbarello, M. Cornalba, P. Griffiths, and J.
Harris for the original concept and general framework of this paper, and E.
Arbarello and M. Cornalba in particular for many helpful conversations
without which this work could not have been completed. Also he wishes to
thank the referee and R. Miranda for ferreting out an elusive mistake in a
previous version of this paper, E. Casini and C. Hacon for help with the
rewrite (especially for pointing me toward Lemma \ref{1.9}), and the Scuola
Normale Superiore, Pisa, Italia, for its hospitality and support during part
of the period of this research.

\section{Deformations of manifolds and differential operators\label{1}}

\subsection{Review of formal Kuranishi theory}

We begin with a brief review of the Newlander-Nirenberg-Kuranishi theory of
deformations of complex structures (see \cite{Ku}, \cite{Ko}, II.1 of \cite
{Gr}, or \cite{C2}). Let 
\begin{equation}
M\overset{\pi }{\longrightarrow }\Delta =\left\{ t\in \Bbb{C}:\left|
t\right| <1\right\}  \label{1.1}
\end{equation}
be a deformation of a complex manifold $M_{0}$ of dimension $m$. Since we
are doing formal deformation theory, all calculations will actually take
place over the formal neighborhood of $0$ in $\Delta $. However, convergence
will not be an issue in anything that we do since we will always be working
from a situation in which we are \textit{given} a geometric deformation and
deriving consequences in the category of formal deformations.

\begin{definition}
\label{1.2}A $C^{\infty }$-diffeomorphism 
\begin{equation*}
F_{\sigma }=\left( \sigma ,\pi \right) :M\rightarrow M_{0}\times \Delta 
\end{equation*}
will be called a trivialization of the deformation $M/\Delta $ if 
\begin{equation*}
\left. \sigma \right| _{M_{0}}=identity
\end{equation*}
and 
\begin{equation*}
\sigma ^{-1}\left( x_{0}\right) 
\end{equation*}
is an analytic disk for each $x_{0}\in M_{0}$.
\end{definition}

The next four lemmas are standard from formal Kuranishi theory:

\begin{lemma}
\label{1.3}Let 
\begin{equation*}
T_{M_{0}}^{*}
\end{equation*}
denote the complexification of the real cotangent bundle of $M_{0}$. Given
any trivialization $F_{\sigma }$, the holomophic cotangent bundle of $M_{t}$
under the $C^{\infty }$-isomorphisms 
\begin{equation*}
M_{t}\cong M_{0}
\end{equation*}
induced by $F_{\sigma }$ corresponds to a subbundle 
\begin{equation*}
T_{t}^{1,0}\subseteq T_{M_{0}}^{*}.
\end{equation*}
If 
\begin{equation*}
\pi ^{1,0}+\pi ^{0,1}:T_{M_{0}}^{*}\rightarrow T_{M_{0}}^{1,0}\oplus
T_{M_{0}}^{0,1}
\end{equation*}
are the two projections, the retriction 
\begin{equation*}
\pi ^{1,0}:T_{t}^{1,0}\rightarrow T_{M_{0}}^{1,0}
\end{equation*}
is an isomorphism for small $t$ so that the composition 
\begin{equation*}
T_{M_{0}}^{1,0}\overset{\left( \pi ^{1,0}\right) ^{-1}}{\longrightarrow }%
T_{t}^{1,0}\overset{\pi ^{0,1}}{\longrightarrow }T_{M_{0}}^{0,1},
\end{equation*}
gives a $C^{\infty }$-mapping 
\begin{equation*}
\xi \left( t\right) :T_{M_{0}}^{1,0}\rightarrow T_{M_{0}}^{0,1}
\end{equation*}
which determines the deformation of (almost) complex structure.
\end{lemma}

Thus, at least formally, we can write 
\begin{equation*}
\xi \left( t\right) =\sum\nolimits_{i>0}\xi _{i}t^{i}
\end{equation*}
with each $\xi _{i}\in A_{M_{0}}^{0,1}\left( T_{1,0}\right) ,$ that is, each 
$\xi _{i}$ is a $\left( 0,1\right) $-form with coefficients in the
holomorphic tangent bundle $T_{1,0}$ of $M_{0}$.

\begin{lemma}
\label{1.4}Every relative complex-valued $C^{\infty }$ -differential form $%
\omega $ on $M/\Delta $ of type $\left( 0,q\right) $ corresponds on a
(formal) neighborhood of $M_{0}$ to a form 
\begin{equation*}
\pi ^{0,q}\left( \left( F_{\sigma }^{-1}\right) ^{*}\left( \omega \right)
\right) =\sum\nolimits_{i,j=0}^{\infty }\omega _{i,j}t^{i}\overline{t}^{j}
\end{equation*}
on 
\begin{equation*}
M_{0}\times \Delta 
\end{equation*}
and so, working modulo $\overline{t}$, gives a holomorphic family 
\begin{equation*}
\omega _{\sigma }:=\sum\nolimits_{i=0}^{\infty }\omega _{i,0}t^{i}.
\end{equation*}
of $C^{\infty }$-forms. This correspondence is a formal isomorphism 
\begin{equation*}
\left( \ \right) _{\sigma }^{q}:\frac{A_{M/\Delta }^{0,q}}{\left\{ \overline{%
t}\right\} }\rightarrow A_{M_{0}}^{0,q}\otimes \Bbb{C}\left[ \left[ t\right]
\right] .
\end{equation*}
If we have two different trivializations $\sigma $ and $\sigma ^{\prime }$,
we have a formal isomorphism
\end{lemma}

\begin{equation*}
G_{\sigma ^{\prime }\sigma }^{q}=\left( \ \right) _{\sigma ^{\prime
}}^{q}\circ \left( \left( \ \right) _{\sigma }^{q}\right) ^{-1}.
\end{equation*}

\begin{lemma}
\label{1.5}For any $C^{\infty }$-function $f$ on $M$ write 
\begin{equation*}
f\circ F_{\sigma }^{-1}=\sum\nolimits_{i,j=0}^{\infty }f_{i,j}t^{i}\overline{%
t}^{j}
\end{equation*}
and define as above 
\begin{equation*}
f_{\sigma }=\sum\nolimits_{i,j=0}^{\infty }f_{i,0}t^{i}.
\end{equation*}
Then define 
\begin{equation*}
\bar{D}_{\sigma }\left( f_{\sigma }\right) :=\left( \bar{\partial }%
_{M_{0}}-\sum\nolimits_{j=1}^{\infty }t^{j}\xi _{j}\right) \left( f_{\sigma
}\right) =\sum\nolimits_{i=0}^{\infty }\bar{\partial }_{M_{0}}f_{i,\sigma
}t^{i}-\sum\nolimits_{i=0,j=1}^{\infty ,\infty }\xi _{j}\left( f_{i,\sigma
}\right) t^{i+j}.
\end{equation*}
Then 
\begin{equation*}
\left( \ \right) _{\sigma }^{1}\circ \bar{\partial }_{M}=\bar{D}_{\sigma
}\circ \left( \ \right) _{\sigma }^{0},
\end{equation*}
and 
\begin{equation*}
\bar{D}_{\sigma }=G_{\sigma \sigma ^{\prime }}^{1}\circ \bar{D}_{\sigma
^{\prime }}\circ G_{\sigma ^{\prime }\sigma }^{0}.
\end{equation*}
Also 
\begin{equation*}
f_{\sigma }\circ F_{\sigma }
\end{equation*}
is holomorphic on $M$ if and only if 
\begin{equation*}
\bar{D}_{\sigma }\left( f_{\sigma }\right) =0.
\end{equation*}
\end{lemma}

We next ask which sequences $\xi _{j}\in A^{0,1}\left( T_{1,0}\right) $ come
from a trivialization of a deformation $\left( \ref{1.1}\right) $. Before
answering this question, we need to make precise the various actions of an
element $\xi \in A^{0,k}\left( T_{1,0}\right) $ on $\sum A^{p,q}\left(
M_{0}\right) $. For any we write the action via contraction as 
\begin{equation*}
\left\langle \xi \right. \left| \ \right\rangle ,
\end{equation*}
and ``Lie differentiation'' as 
\begin{equation*}
L_{\xi }:=\left\langle \xi \right. \left| \ \right\rangle \circ d+\left(
-1\right) ^{k}d\circ \left\langle \xi \right. \left| \ \right\rangle .
\end{equation*}
The sign is so chosen that, writing any element of $A^{0,k}\left(
T_{1,0}\right) $ locally as a sum of terms 
\begin{equation*}
\xi =\bar{\eta}\otimes \chi
\end{equation*}
for some closed $\left( 0,k\right) $-form $\bar{\eta}$ and $\chi \in
A^{0,0}\left( T_{1,0}\right) $, then 
\begin{equation*}
L_{\xi }=\bar{\eta}\otimes L_{\chi }.
\end{equation*}
(Warning: Since, as an operator on $A^{0,q}\left( M_{0}\right) $, $L_{f\xi
}=fL_{\xi }$, one has 
\begin{equation*}
\left[ \overline{\partial },L_{\xi }\right] =L_{\overline{\partial }\xi
}:A^{0,q}\left( M_{0}\right) \rightarrow A^{0,q+k+1}\left( M_{0}\right)
\end{equation*}
however the identity does not hold as an operator on $A^{p,q}\left(
M_{0}\right) $ for $p>0$.)

Also we compute 
\begin{eqnarray*}
&&L_{\xi }L_{\xi ^{\prime }}-\left( -1\right) ^{\deg \bar{\eta}\cdot \deg 
\bar{\eta}^{\prime }}L_{\xi ^{\prime }}L_{\xi } \\
&=&\left( \bar{\eta}\otimes L_{\chi }\right) \left( \bar{\eta}^{\prime
}\otimes L_{\chi ^{\prime }}\right) -\left( -1\right) ^{\deg \bar{\eta}\cdot
\deg \bar{\eta}^{\prime }}\left( \bar{\eta}^{\prime }\otimes L_{\chi
^{\prime }}\right) \left( \bar{\eta}\otimes L_{\chi }\right) \\
&=&\bar{\eta}\bar{\eta}^{\prime }\left( L_{\chi }L_{\chi ^{\prime }}-L_{\chi
^{\prime }}L_{\chi }\right) \\
&=&\bar{\eta}\bar{\eta}^{\prime }L_{\left[ \chi ,\chi ^{\prime }\right] }.
\end{eqnarray*}
So, using this local presentation for 
\begin{equation*}
\xi \in A^{0,j}\left( T_{1,0}\right) ,\xi ^{\prime }\in A^{0,k}\left(
T_{1,0}\right) ,
\end{equation*}
we can define 
\begin{equation*}
\left[ \xi ,\xi ^{\prime }\right] =\bar{\eta}\bar{\eta}^{\prime }\left[ \chi
,\chi ^{\prime }\right] \in A^{0,j+k}\left( T_{1,0}\right) .
\end{equation*}

\begin{lemma}
\label{1.6}The almost complex structures given on a coordinate neighborhood $%
W_{0}$ in $M_{0}$ by the the $\left( 0,1\right) $-tangent distributions 
\begin{equation*}
\left( \frac{\partial }{\partial \overline{v_{W_{0}}^{k}}}%
-\sum\nolimits_{i=1}^{\infty }\sum\nolimits_{l}h_{i,k}^{l}t^{i}\frac{%
\partial }{\partial v_{W_{0}}^{l}}\right) 
\end{equation*}
are integrable, that is, come from a deformation/trivialization of $M_{0}$
as in Definition $\left( \ref{1.2}\right) $, if and only if, for 
\begin{equation*}
\xi =\sum\nolimits_{i=1}^{\infty }\sum\nolimits_{k,l}d\overline{v_{W_{0}}^{k}%
}\wedge h_{i,k}^{l}t^{i}\frac{\partial }{\partial v_{W_{0}}^{l}},
\end{equation*}
we have 
\begin{equation*}
\bar{\partial }\xi =\frac{1}{2}\left[ \xi ,\xi \right] .
\end{equation*}
\end{lemma}

\begin{proposition}
\label{1.7}Two trivializations $F_{\sigma }$ and $F_{\sigma ^{\prime }}$ of
the same deformation $\left( \ref{1.1}\right) $ are related by a holomorphic
automorphism $\varphi $ of $M/\Delta ,$ that is, there is a commutative
diagram 
\begin{equation*}
\begin{tabular}{ccc}
$M$ & $\overset{\varphi }{\longrightarrow }$ & $M$ \\ 
$\downarrow \sigma $ &  & $\downarrow \sigma ^{\prime }$ \\ 
$M_{0}$ & $=$ & $M_{0}$%
\end{tabular}
,
\end{equation*}
if and only if 
\begin{equation*}
\bar{D}_{\sigma }=\bar{D}_{\sigma ^{\prime }}.
\end{equation*}
\end{proposition}

\begin{proof}
One implication is immediate from the definitions of $\bar{D}_{\sigma }$ and 
$\bar{D}_{\sigma ^{\prime }}$. For the other, the equality 
\begin{equation*}
\xi _{\sigma }=\xi _{\sigma ^{\prime }}
\end{equation*}
implies that the differential of the $C^{\infty }$-automorphism 
\begin{equation*}
\varphi :=\left( \sigma ,\pi \right) ^{-1}\circ \left( \sigma ,\pi \right)
:M\rightarrow M
\end{equation*}
preserves the $\left( 1,0\right) $-subspace of the (complexified) tangent
space and therefore $\varphi $ is holomorphic.
\end{proof}

\subsection{Gauge transformations\label{1.8}}

We begin now with a deformation 
\begin{equation*}
M/\Delta
\end{equation*}
of $M_{0}$ and let 
\begin{equation*}
F_{\sigma }:M\overset{\left( \sigma ,\pi \right) }{\longrightarrow }%
M_{0}\times \Delta
\end{equation*}
be a trivialization with associated Kuranishi data 
\begin{equation*}
\xi _{\sigma }.
\end{equation*}
Suppose that we have a one-real-parameter group of diffeomorphisms 
\begin{equation*}
\Phi _{s}:M_{0}\times \Delta \rightarrow M_{0}\times \Delta
\end{equation*}
defined over $\Delta $ such that 
\begin{equation*}
F_{s}:=\Phi _{s}\circ F_{\sigma }:M\overset{\left( \sigma _{s},\pi \right) }{%
\longrightarrow }M_{0}\times \Delta
\end{equation*}
is a trivialization for each (sufficiently small) $s$ and, for each $%
x_{0}\in M_{0}$, 
\begin{equation*}
\left. \Phi _{s}\right| _{\left\{ x_{0}\right\} \times \Delta }
\end{equation*}
is a real-analytic family of complex-analytic embeddings of $\Delta $ in $%
M_{0}\times \Delta $. Then, as for example in \S 2 of \cite{C2}, there is
then associated a vector field 
\begin{equation*}
\beta +\overline{\beta }
\end{equation*}
where 
\begin{equation*}
\beta =\sum\nolimits_{j>0}\beta _{j}t^{j}
\end{equation*}
and each $\beta _{j}$ is a $C^{\infty }$-vector field of type $\left(
1,0\right) $ on $M_{0}$, such that, for 
\begin{equation*}
g=g_{0}+g_{1}t+\ldots ,
\end{equation*}
on $M_{0}\times \Delta $ we have 
\begin{equation}
g\circ \Phi _{s}=e^{L_{s\beta +\overline{s\beta }}}\left( g\right) .
\label{1.8.4}
\end{equation}
We let 
\begin{equation*}
F_{\beta }:=F_{1}=\Phi _{1}\circ F_{\sigma }:M\rightarrow M_{0}\times \Delta
.
\end{equation*}
Then.by $\left( \ref{1.8.4}\right) $ we have for any $C^{\infty }$-function $%
g$ on $M$ that 
\begin{equation}
g_{\beta }=e^{L_{\beta }}\left( g_{\sigma }\right) .  \label{1.8.3}
\end{equation}

If $\xi _{s}$ denotes the Kuranishi data for the trivialization $F_{s}$,
then by direct computation 
\begin{equation*}
\frac{\partial \xi _{s}}{\partial s}=\left[ \overline{\partial },\beta
\right] +\left[ \beta ,\xi _{s}\right] .
\end{equation*}
(See for example \S 3 of \cite{C2}.) On the other hand, if we define 
\begin{equation}
\varsigma _{\beta }:=\frac{e^{\left[ \beta ,\ \right] }-1}{\left[ \beta ,\
\right] }\left( \left[ \overline{\partial },\beta \right] \right)
\label{1.8.1}
\end{equation}
and the action 
\begin{equation}
\xi _{\beta }:=\beta \cdot \left( \xi \right) =e^{\left[ \beta ,\ \right]
}\left( \xi \right) +\varsigma _{\beta }  \label{1.8.5}
\end{equation}
one also has by direct computation that 
\begin{equation*}
\frac{\partial \xi _{s\beta }}{\partial s}=\left[ \overline{\partial },\beta
\right] +\left[ \beta ,\xi _{s\beta }\right] .
\end{equation*}
(See for example \S 3 of \cite{C2}. Compare with \S 3 of \cite{GM}.) The
conclusion is that $\xi _{s\beta }$ is the Kuranishi data for the
trivialization $F_{s}$ for all $s$ and so, in particular 
\begin{equation*}
\xi _{\beta }
\end{equation*}
is the Kuranishi data for the trivialization $F_{1}=F_{\beta }$.

So the group of vector fields $\beta $ acts on the Kuranishi data associated
to the deformation $M/\Delta $. This action corresponds to the change of the
given trivialization by a $C^{\infty }$-automorphism 
\begin{equation}
\Phi _{\beta }:M_{0}\times \Delta \rightarrow M_{0}\times \Delta
\label{1.8.2}
\end{equation}
defined over $\Delta $.

\begin{lemma}
\label{1.9}i) 
\begin{equation*}
\left( e^{L_{\beta }}\right) \left( \overline{\partial }-\xi \right) \left(
e^{-L_{\beta }}\right) =\overline{\partial }-\xi _{\beta }.
\end{equation*}
ii) Given a function 
\begin{equation*}
f_{\beta }=\sum\nolimits_{i}f_{\beta ,i}t^{i}
\end{equation*}
on $M_{0}\times \Delta ,$ the function 
\begin{equation*}
f_{\beta }\circ F_{\beta }
\end{equation*}
is holomorphic on $M$ if and only if 
\begin{equation*}
\left( \overline{\partial }-\xi \right) \left( e^{-L_{\beta }}\left(
f_{\beta }\right) \right) =0.
\end{equation*}
\end{lemma}

\begin{proof}
i) This assertion is implicit in $\left( \ref{1.8.3}\right) $ but, as a
check, we will do it by direct comptation. 
\begin{equation*}
\left( e^{L_{\beta }}\right) \left( \overline{\partial }-\xi \right) \left(
e^{-L_{\beta }}\right) =\overline{\partial }+\left[ e^{L_{\beta }},\overline{
\partial }\right] \left( e^{-L_{\beta }}\right) -e^{\left[ \beta ,\ \right]
}\left( \xi \right) .
\end{equation*}
If we can show the identity 
\begin{equation}
\left[ \overline{\partial },e^{L_{\beta }}\right] =\varsigma _{\beta }\circ
e^{L_{\beta }},  \label{1.9.1}
\end{equation}
the lemma will follow from Lemma \ref{1.5} since, by definition, 
\begin{equation*}
\xi _{\beta }=e^{\left[ \beta ,\ \right] }\left( \xi \right) +\varsigma
_{\beta }.
\end{equation*}
To see $\left( \ref{1.9.1}\right) $ we prove by induction that 
\begin{equation*}
\left[ \overline{\partial },\beta ^{n+1}\right]
=\sum\nolimits_{i=0}^{n}\left( 
\begin{array}{c}
n+1 \\ 
i
\end{array}
\right) \left( \left[ \beta ,\ \right] \right) ^{n-i}\left[ \overline{%
\partial },\beta \right] \beta ^{i}.
\end{equation*}
Inductively 
\begin{eqnarray*}
\left[ \overline{\partial },\beta ^{n+1}\right] &=&\left[ \overline{\partial 
},\beta \right] \cdot \beta ^{n}+\beta \cdot \left[ \overline{\partial }
,\beta ^{n}\right] \\
&=&\left[ \overline{\partial },\beta \right] \cdot \beta ^{n}+\beta \left(
\left( 
\begin{array}{c}
n \\ 
0
\end{array}
\right) \left[ \beta ,\ \right] ^{n-1}\cdot \left[ \overline{\partial }%
,\beta \right] +\ldots +\left( 
\begin{array}{c}
n \\ 
n-1
\end{array}
\right) \left[ \overline{\partial },\beta \right] \cdot \beta ^{n-1}\right)
\\
&=&\left( \left( 
\begin{array}{c}
n \\ 
n
\end{array}
\right) \left[ \overline{\partial },\beta \right] \cdot \beta ^{n}+\left( 
\begin{array}{c}
\left( 
\begin{array}{c}
n \\ 
0
\end{array}
\right) \left( \left( \left[ \beta ,\ \right] \right) ^{n}\left[ \overline{%
\partial },\beta \right] +\left( \left[ \beta ,\ \right] \right)
^{n-1}\left[ \overline{\partial },\beta \right] \beta \right) \\ 
+\ldots + \\ 
\left( 
\begin{array}{c}
n \\ 
n-1
\end{array}
\right) \left( \left( \left[ \beta ,\ \right] \right) \left[ \overline{%
\partial },\beta \right] \beta ^{n-1}+\left[ \overline{\partial },\beta
\right] \beta ^{n}\right)
\end{array}
\right) \right)
\end{eqnarray*}
Now use the identity 
\begin{equation*}
\left( \left( 
\begin{array}{c}
n-1 \\ 
r
\end{array}
\right) +\left( 
\begin{array}{c}
n-1 \\ 
r-1
\end{array}
\right) \right) =\left( 
\begin{array}{c}
n \\ 
r
\end{array}
\right) .
\end{equation*}
Thus 
\begin{eqnarray*}
\left[ \overline{\partial },\left( \sum\nolimits_{n=0}^{\infty }\frac{\beta
^{n}}{n!}\right) \right] &=&\sum\nolimits_{n=1}^{\infty }\frac{1}{n!}%
\sum\nolimits_{i=0}^{n-1}\left( 
\begin{array}{c}
n \\ 
i
\end{array}
\right) \left( \left[ \beta ,\ \right] \right) ^{n-1-i}\left[ \overline{%
\partial },\beta \right] \beta ^{i} \\
&=&\sum\nolimits_{k=0,i=0}^{\infty ,\infty }\frac{1}{\left( k+1\right) !}%
\left( \left[ \beta ,\ \right] \right) ^{k}\left[ \overline{\partial },\beta
\right] \frac{\beta ^{i}}{i!} \\
&=&\xi _{\beta }\circ \left( \sum\nolimits_{i=0}^{\infty }\frac{\beta ^{i}}{%
i!}\right) .
\end{eqnarray*}

ii) 
\begin{eqnarray*}
0 &=&\left[ \overline{\partial },e^{L_{\beta }}\circ e^{-L_{\beta }}\right]
\\
&=&\left[ \overline{\partial },e^{L_{\beta }}\right] \circ e^{-L_{\beta
}}+e^{L_{\beta }}\circ \left[ \overline{\partial },e^{-L_{\beta }}\right] \\
&=&\varsigma _{\beta }+e^{L_{\beta }}\circ \varsigma _{-\beta }\circ
e^{-L_{\beta }}
\end{eqnarray*}
so that 
\begin{equation*}
e^{L_{\beta }}\circ \left( \varsigma _{-\beta }-\xi \right) \circ
e^{-L_{\beta }}=-\left( \varsigma _{\beta }+\left( e^{L_{\beta }}\circ \xi
\circ e^{-L_{\beta }}\right) \right) .
\end{equation*}
\end{proof}

Suppose now that we have two trivializations 
\begin{eqnarray*}
F_{\sigma } &:&M\rightarrow M_{0}\times \Delta \\
F_{\sigma ^{\prime }} &:&M\rightarrow M_{0}\times \Delta
\end{eqnarray*}
of a given deformation 
\begin{equation*}
M/\Delta .
\end{equation*}
Then 
\begin{equation*}
F_{\sigma ^{\prime }}\circ F_{\sigma }^{-1}
\end{equation*}
is a $C^{\infty }$-diffeomorphism of $M_{0}\times \Delta $ and so can be
realized as the value at $s=1$ of a one-parameter group of diffeomorphisms
which restrict to an analytic family of analytic embbeddings of $\left\{
x_{0}\right\} \times \Delta $ for each $x_{0}\in M_{0}$. Thus referring to
the notation of Lemma \ref{1.5} above we have that there is a $C^{\infty }$%
-vector field $\kappa $ of type $\left( 1,0\right) $ such that 
\begin{equation}
\begin{array}{l}
g_{\sigma }=e^{L_{-\kappa }}\left( g_{\sigma ^{\prime }}\right) \\ 
\overline{D}_{\sigma ^{\prime }}=e^{L_{\kappa }}\circ \overline{D}_{\sigma
}\circ e^{L_{-\kappa }}.
\end{array}
\label{1.10}
\end{equation}

\subsection{Schiffer-type deformations\label{1.11}}

We now consider a special class of deformations of $M_{0}$, those for which
the change of complex structure can be localized at a union $A_{0}$ of
codimension-one subvarieties on $M_{0}$. We let 
\begin{equation}
\beta \in A_{M_{0}}^{0,0}\left( T_{M_{0}}^{1,0}\right) \otimes t\Bbb{C}%
\left[ \left[ t\right] \right]  \label{1.11.1}
\end{equation}
be a vector field which is

i) meromorphic in an analytic neighborhood $\left( U_{0}\times \Delta
\right) $ of the set $\left( A_{0}\times \Delta \right) $ on $\left(
M_{0}\times \Delta \right) $,

ii) $C^{\infty }$ on $\left( M_{0}-A_{0}\right) \times \Delta $.

\noindent Using Lemma \ref{1.9} for the case in which we first take 
\begin{equation}
F_{\sigma }:\left( \left( M_{0}-A_{0}\right) \times \Delta \right)
\rightarrow \left( \left( M_{0}-A_{0}\right) \times \Delta \right)
\label{1.11.4}
\end{equation}
in \ref{1.8} as the identity map, we define a deformation $M_{\beta }/\Delta 
$ of $M_{0}$ by the integrable Kuranishi data 
\begin{equation}
\xi _{\beta }:=\varsigma _{\beta }.  \label{1.11.2}
\end{equation}
Notice that $\xi _{\beta }=0$ in a neighborhood of $A_{0}\times \Delta $ so $%
\xi _{\beta }$ corresponds to a trivialization 
\begin{equation*}
F_{\beta }:M_{\beta }\overset{\left( \sigma _{\beta },\pi \right) }{%
\longrightarrow }M_{0}\times \Delta
\end{equation*}
with 
\begin{equation*}
F_{\beta }:\left( \sigma _{\beta }\right) ^{-1}\left( U_{0}\right)
\rightarrow U_{0}\times \Delta
\end{equation*}
an analytic isomorphism. Denote 
\begin{equation}
\overline{D}_{\beta }:=\overline{D}_{\sigma _{\beta }}=\overline{\partial }%
-\varsigma _{\beta }.  \label{1.11.3}
\end{equation}
We call $A_{0}$ the \textit{center} of the Schiffer-type deformation.

Let 
\begin{equation*}
A_{\beta }:=\left( \sigma _{\beta }\right) ^{-1}\left( A_{0}\right)
\subseteq M_{\beta }.
\end{equation*}
From \ref{1.8}, Lemma \ref{1.9} and the above we conclude:

\begin{lemma}
\label{1.12} 
\begin{equation*}
f_{\beta }\circ F_{\beta }
\end{equation*}
is analytic on $M_{\beta }$ if and only if 
\begin{equation*}
\overline{\partial }_{M_{0}}\left( e^{-L_{\beta }}\left( f_{\beta }\right)
\right) =0.
\end{equation*}
\end{lemma}

In fact, for any divisor $B_{0}$ supported on $A_{0}$, $B_{0}$ has a unique
extension to a divisor 
\begin{equation*}
B_{\beta }
\end{equation*}
on $M_{\beta }$ which is supported on $A_{\beta }$. We denote by 
\begin{equation*}
C_{B_{0}}
\end{equation*}
the vector space of functions $f_{0}$ which are $C^{\infty }$ on $\left(
M_{0}-A_{0}\right) $ and meromorphic on $U_{0}$ and for which 
\begin{equation*}
B_{0}+div\left( f_{0}\right)
\end{equation*}
is effective on $U_{0}$. Then:

\begin{lemma}
\label{1.13}i) A meromorphic function $f$ on $M_{\beta }$ with 
\begin{equation*}
B_{\beta }+div\left( f\right) 
\end{equation*}
effective is a formal sum 
\begin{equation*}
f_{\beta }:=f_{\sigma _{\beta }}=\sum\nolimits_{i=0}^{\infty }f_{\beta
,i}t^{i}
\end{equation*}
such that each $f_{\beta ,i}\in C_{B_{0}}$ and 
\begin{equation*}
\left( \overline{\partial }-\sum\nolimits_{j=1}^{\infty }\xi _{\beta
,j}t^{j}\right) \left( \sum\nolimits_{i=0}^{\infty }f_{\beta ,i}t^{i}\right)
=0.
\end{equation*}
ii) The meromorphic functions on $M_{\beta }$ with 
\begin{equation*}
B_{\beta }+div\left( f\right) 
\end{equation*}
effective are given by the kernel of the mapping 
\begin{equation*}
e^{L_{\beta }}:H^{0}\left( \mathcal{O}_{M_{0}}\left( \infty \cdot %
A_{0}\right) \right) \otimes \Bbb{C}\left[ \left[ t\right] \right]
\rightarrow H^{0}\left( \frac{\mathcal{O}_{M_{0}}\left( \infty \cdot %
A_{0}\right) }{\mathcal{O}_{M_{0}}\left( B_{0}\right) }\right) \otimes \Bbb{C%
}\left[ \left[ t\right] \right] .
\end{equation*}
iii) If 
\begin{equation*}
i:A_{0}\rightarrow M_{0}
\end{equation*}
is the inclusion map and $R$ denote the image of the map 
\begin{equation*}
\left( i^{-1}\mathcal{O}_{M_{0}}\left( B_{0}\right) \otimes t\Bbb{C}\left[
\left[ t\right] \right] \right) \overset{\overline{\partial }\circ
e^{-L_{\beta }}}{\longrightarrow }\left( H^{1}\left( \mathcal{O}%
_{M_{0}}\left( B_{0}\right) \right) \otimes \Bbb{C}\left[ \left[ t\right]
\right] \right) ,
\end{equation*}
then $f_{0}\in H^{0}\left( \mathcal{O}_{M_{0}}\left( B_{0}\right) \right) $
extends to a global section of $\mathcal{O}_{M_{\beta }}\left( B_{\beta
}\right) $ if and only if 
\begin{equation*}
\left[ \overline{\partial },e^{-L_{\beta }}\right] \left( f_{0}\right) \in R.
\end{equation*}
\end{lemma}

\begin{proof}
i) The assertion is immediate from Lemma \ref{1.5}.

ii) Again by Lemma \ref{1.9}i) occurs exactly when $f_{\beta }$ lies in 
\begin{equation*}
C_{B_{0}}\otimes \Bbb{C}\left[ \left[ t\right] \right] \cap image\left(
H^{0}\left( \mathcal{O}_{M_{0}}\left( \infty \cdot A_{0}\right) \right)
\otimes \Bbb{C}\left[ \left[ t\right] \right] \overset{e^{L_{\beta }}}{%
\longrightarrow }C_{\infty \cdot A_{0}}\otimes \Bbb{C}\left[ \left[ t\right]
\right] \right) .
\end{equation*}

iii) follows Lemma \ref{1.12} and from the cohomology exact sequence
associated to the short exact sequence 
\begin{equation*}
0\rightarrow \mathcal{O}_{M_{0}}\left( B_{0}\right) \rightarrow \mathcal{O}
_{M_{0}}\left( \infty \cdot A_{0}\right) \rightarrow \frac{\mathcal{O}%
_{M_{0}}\left( \infty \cdot A_{0}\right) }{\mathcal{O}_{M_{0}}\left(
B_{0}\right) }\rightarrow 0.
\end{equation*}
\end{proof}

\subsection{Gauge transformation on Schiffer-type trivializations\label{1.14}
}

Next suppose we wish to change our trivialization 
\begin{equation*}
F_{\beta }:M_{\beta }\overset{\left( \sigma _{\beta },\pi \right) }{%
\longrightarrow }M_{0}\times \Delta
\end{equation*}
by an allowable $C^{\infty }$-automorphism 
\begin{equation*}
\begin{array}{ccc}
M_{\beta } & \overset{F_{\beta }}{\longrightarrow } & M_{0}\times \Delta \\ 
= &  & \downarrow \Phi \\ 
M_{\beta } & \overset{G}{\longrightarrow } & M_{0}\times \Delta
\end{array}
\end{equation*}
defined over $\Delta $. That is

\begin{enumerate}
\item  such that $\Phi $ preserves $A_{0}\times \Delta $ as a set

\item  is holomorphic on $U_{0}\times \Delta $.

\item  $\Phi $ restricts to an analytic embedding of each disk $\left\{
x_{0}\right\} \times \Delta $.
\end{enumerate}

To calculate the Kuranishi data for $G,$ we proceed as in \ref{1.8}. We can
assume that $\Phi =\Phi _{1}$ for a family $\Phi _{s}$ as in \ref{1.8}. We
can further assume that $\left. \Phi _{s}\right| _{U_{0}\times \Delta }$ is
a real analytic family of complex analytic maps. Let $\kappa
=\sum\nolimits_{j=1}^{\infty }\kappa _{j}t^{j}$ denote the $C^{\infty }$
-vector field of type $\left( 1,0\right) $ such that the family $\Phi _{s}$
is associated to 
\begin{equation*}
s\left( \kappa +\overline{\kappa }\right) .
\end{equation*}
Then by $\left( \ref{1.10}\right) $ for $F_{\sigma }=F_{\beta }$ and $%
F_{\sigma ^{\prime }}=G$ we have 
\begin{equation*}
\begin{array}{l}
g_{\beta }=e^{L_{-\kappa }}\left( g_{\sigma ^{\prime }}\right) \\ 
\overline{D}_{\sigma ^{\prime }}=e^{L_{\kappa }}\circ \overline{D}_{\beta
}\circ e^{L_{-\kappa }}.
\end{array}
\end{equation*}

Computing using $\left( \ref{1.8.5}\right) $ and $\left( \ref{1.9.1}\right) $%
\begin{eqnarray*}
e^{L_{\kappa }}\circ \overline{D}_{\beta }\circ e^{L_{-\kappa }}
&=&e^{L_{\kappa }}\circ \left( \overline{\partial }-\varsigma _{\beta
}\right) \circ e^{L_{-\kappa }} \\
&=&e^{L_{\kappa }}\circ \left( \overline{\partial }-\left[ \overline{%
\partial },e^{L_{\beta }}\right] \circ e^{L_{-\beta }}\right) \circ
e^{L_{-\kappa }} \\
&=&e^{L_{\kappa }}\circ \left( e^{L_{\beta }}\circ \overline{\partial }\circ
e^{L_{-\beta }}\right) \circ e^{L_{-\kappa }}.
\end{eqnarray*}
Thus we conclude that $g_{\sigma ^{\prime }}$ is holomorphic if and only if 
\begin{equation*}
\overline{\partial }_{M_{0}}\left( e^{L_{-\beta }}\circ e^{L_{-\kappa
}}\left( g_{\sigma ^{\prime }}\right) \right) =0.
\end{equation*}

\begin{lemma}
\label{1.15}For a power series 
\begin{equation*}
g=\sum\nolimits_{i=0}^{\infty }g_{i}t^{i}
\end{equation*}
on $M_{0}\times \Delta $, $g\circ G$ is holomorphic on $M_{\beta }$ if and
only if 
\begin{equation*}
\overline{\partial }_{M_{0}}\left( e^{L_{-\beta }}\circ e^{L_{-\kappa
}}\left( g\right) \right) =0.
\end{equation*}
\end{lemma}

\section{Deformations of line bundles and differential operators\label{2}}

\subsection{The $\mu $-maps\label{2.1}}

Let $X_{0}$ be a complex manifold and let $L_{0}$ be a holomorphic line
bundle on $X_{0}$. Let 
\begin{equation*}
\frak{D}\left( L_{0}\right) ,\frak{D}_{n}\left( L_{0}\right)
\end{equation*}
denote the sheaf of (holomorphic) differential operators, respectively the
sheaf of differential operators of order $\leq n$, on (sections of) the line
bundle $L_{0}$. Whenever 
\begin{equation*}
H^{2}\left( \frak{D}_{n}\left( L_{0}\right) \right) =0
\end{equation*}
we have a natural exact sequence 
\begin{equation*}
H^{1}\left( \frak{D}_{n}\left( L_{0}\right) \right) \rightarrow H^{1}\left( 
\frak{D}_{n+1}\left( L_{0}\right) \right) \rightarrow H^{1}\left(
S^{n+1}T_{X_{0}}\right) \rightarrow 0
\end{equation*}
where the second last map is induced by the symbol map on differential
operators. So there exists natural mappings 
\begin{equation}
\tilde{\mu}^{n}:H^{1}\left( \frak{D}_{n}\left( L_{0}\right) \right)
\rightarrow Hom\left( H^{0}\left( L_{0}\right) ,H^{1}\left( L_{0}\right)
\right)  \label{2.1.1}
\end{equation}
and 
\begin{equation}
\mu ^{n+1}:H^{1}\left( S^{n+1}T_{X_{0}}\right) \rightarrow \frac{Hom\left(
H^{0}\left( L_{0}\right) ,H^{1}\left( L_{0}\right) \right) }{image\ \tilde{%
\mu}^{n}}.  \label{2.1.2}
\end{equation}
(In the next chapter we will establish Petri's conjecture on generic curve $%
C_{0}$ by establishing that the mappings $\left( \ref{2.1.2}\right) $ are
zero for $n\geq 0$ and $X_{0}=C_{0}$.)

Suppose now that we are given a deformation 
\begin{equation}
L\overset{p}{\longrightarrow }X\overset{\pi }{\longrightarrow }\Delta
\label{2.2.1}
\end{equation}
of the pair $\left( L_{0},X_{0}\right) $. We consider $C^{\infty }$-sections
of $L$ as $C^{\infty }$-functions on the dual line bundle $L^{\vee }$. These
functions $f$ are characterized by the properties 
\begin{equation}
\begin{tabular}{l}
$\chi \left( f\right) =f$ \\ 
$\bar{\chi}\left( f\right) =0$%
\end{tabular}
\label{2.2.2}
\end{equation}
where $\chi $ is the (holomorphic) Euler vector-field associated with the $%
\Bbb{C}^{*}$-action on $L^{\vee }$.

\subsection{Trivializations of deformations of line bundles\label{2.3}}

We next claim that, given a trivialization $\sigma $ of the deformation $%
X/\Delta $ and given a line bundle $L/X$ we can make compatible
trivializations 
\begin{equation}
\begin{tabular}{ccc}
$L^{\vee }$ & $\overset{F_{\lambda }=\left( \lambda ,\pi \circ q\right) }{%
\longrightarrow }$ & $L_{0}^{\vee }\times \Delta $ \\ 
$\downarrow q$ &  & $\downarrow \left( q_{0},id.\right) $ \\ 
$X$ & $\overset{F_{\sigma }=\left( \sigma ,\pi \right) }{\longrightarrow }$
& $X_{0}\times \Delta $ \\ 
$\downarrow \pi $ &  & $\downarrow $ \\ 
$\Delta $ & $=$ & $\Delta $%
\end{tabular}
\label{2.3.1}
\end{equation}
of the deformation $L^{\vee }/X$ of $L_{0}^{\vee }/X_{0}$ as in Lemma \ref
{1.3} but with the additional property that each fiber of the trivialization
respects the structure of holomorphic line bundles, that is, if we denote by 
$\tau =\tau _{\sigma }$ the lifting of $\frac{\partial }{\partial t}$
induced by the trivialization of $X/\Delta $, then $\tau =\tau _{\lambda }$
for the deformation $L^{\vee }$ of $L_{0}^{\vee }$ is obtained as a lifting
of $\tau _{\sigma }$ such that 
\begin{equation}
\left[ \tau _{\lambda },\chi \right] =0.  \label{2.3.2}
\end{equation}

To see that this is always possible, let $\left\{ W\right\} $ be a covering
of $X$ by coordinate disks and $\left\{ W_{0}\right\} $ the restriction of
this covering to $X_{0}$. We construct a $C^{\infty }$ partition-of-unity $%
\left\{ \rho _{W_{0}}\right\} $ subordinate to the induced covering of $%
X_{0} $. Recall that $L$ is given with respect to the trivialization $\sigma 
$ by holomorphic local patching data 
\begin{eqnarray*}
g^{_{WW^{\prime }}}\left( x\right) &=&\sum g_{i}^{_{WW^{\prime }}}\left(
x_{0}\right) t^{i} \\
&=&g^{W_{0}W_{0}^{\prime }}\left( x_{0}\right) \exp \left(
\sum\nolimits_{j>0}a_{j}^{_{WW^{\prime }}}\left( x_{0}\right) t^{j}\right)
\end{eqnarray*}
where $x_{0}=\sigma \left( x\right) $ and 
\begin{equation*}
\sum\nolimits_{j>0}a_{j}^{_{WW^{\prime }}}\left( x_{0}\right) t^{j}=\log 
\frac{g^{_{WW^{\prime }}}\left( x\right) }{g^{W_{0}W_{0}^{\prime }}\left(
x_{0}\right) }.
\end{equation*}
Notice that, if $V,W,$ and $W^{\prime }$ are three open sets of the cover
which have non-empty intersection, then, for all $j>0$, 
\begin{equation*}
a_{j}^{_{VW}}+a_{j}^{_{WW^{\prime }}}=a_{j}^{_{VW^{\prime }}}.
\end{equation*}
Define the mapping 
\begin{equation*}
L\rightarrow L_{0}
\end{equation*}
over $W_{0}\times \Delta $ by 
\begin{equation}
\left( x,v\right) \mapsto \left( x_{0},\exp \left( \sum\nolimits_{W^{\prime
}}\rho _{W_{0}^{\prime }}\left( x_{0}\right) \left(
\sum\nolimits_{j>0}a_{j}^{_{WW^{\prime }}}\left( x_{0}\right) t^{j}\right)
\right) \cdot v\right) .  \label{2.3.0}
\end{equation}
This map is well defined since, over $V\cap W$ we have 
\begin{equation*}
g^{VW}\left( x\right) =g^{V_{0}W_{0}}\left( x_{0}\right) \exp \left(
\sum\nolimits_{j>0}a_{j}^{_{VW}}\left( x_{0}\right) t^{j}\right)
\end{equation*}
and so 
\begin{eqnarray*}
&&g^{_{VW}}\left( x\right) \cdot \exp \left( \sum\nolimits_{W^{\prime }}\rho
_{W_{0}^{\prime }}\left( x_{0}\right) \left(
\sum\nolimits_{j>0}a_{j}^{_{WW^{\prime }}}\left( x_{0}\right) t^{j}\right)
\right) \\
&=&g^{V_{0}W_{0}}\left( x_{0}\right) \exp \left(
\sum\nolimits_{j>0}a_{j}^{_{VW}}\left( x_{0}\right) t^{j}\right) \cdot \exp
\left( \sum\nolimits_{W^{\prime }}\rho _{W_{0}^{\prime }}\left( x_{0}\right)
\left( \sum\nolimits_{j>0}a_{j}^{_{WW^{\prime }}}\left( x_{0}\right)
t^{j}\right) \right) \\
&=&g^{V_{0}W_{0}}\left( x_{0}\right) \exp \left( \sum\nolimits_{W^{\prime
}}\rho _{W_{0}^{\prime }}\left( x_{0}\right) \sum\nolimits_{j>0}\left(
a_{j}^{_{VW}}+a_{j}^{_{WW^{\prime }}}\right) \left( x_{0}\right) t^{j}\right)
\\
&=&g^{V_{0}W_{0}}\left( x_{0}\right) \exp \left( \sum\nolimits_{W^{\prime
}}\rho _{W_{0}^{\prime }}\left( x_{0}\right)
\sum\nolimits_{j>0}a_{j}^{_{VW^{\prime }}}\left( x_{0}\right) t^{j}\right) .
\end{eqnarray*}

Referring to Lemma \ref{1.5} our deformation/trivialization $\left( \ref
{2.3.1}\right) $ is given by 
\begin{equation*}
\xi _{j}\in A^{0,1}\left( T_{L_{0}^{\vee }}\right)
\end{equation*}
for which 
\begin{equation}
L_{\chi }\xi _{j}=L_{\bar{\chi}}\xi _{j}=0.  \label{2.3.3}
\end{equation}

We call a trivialization satisfying $\left( \ref{2.3.1}\right) $-$\left( \ref
{2.3.3}\right) $ a\textit{\ trivialization of line bundles}. We say that the
trivializations $\lambda $ of $L^{\vee }/\Delta $ and $\sigma $ of $X/\Delta 
$ are \textit{compatible} if they make the diagram $\left( \ref{2.3.1}%
\right) $ commutative. By an elementary computation in local coordinates,
sections 
\begin{equation*}
\xi _{i}\in A_{L_{0}^{\vee }}^{0,1}\otimes T_{L_{0}^{\vee }}
\end{equation*}
associated to a trivialization of line bundles lie in a subspace 
\begin{equation*}
A\subseteq A_{L_{0}^{\vee }}^{0,1}\otimes T_{L_{0}^{\vee }}
\end{equation*}
comprising the the middle term of an exact sequence 
\begin{equation}
0\rightarrow q_{0}^{-1}\left( A_{X_{0}}^{0,1}\right) \otimes _{\Bbb{C}}\Bbb{C%
}\chi \rightarrow A\rightarrow q_{0}^{-1}\left( A_{X_{0}}^{0,1}\otimes
T_{X_{0}}\right) \rightarrow 0,  \label{2.3.4}
\end{equation}
that is, 
\begin{equation}
A=A_{X_{0}}^{0,1}\left( \frak{D}_{1}\left( L_{0}\right) \right) .
\label{2.3.5}
\end{equation}

Notice that the first form 
\begin{equation*}
\xi _{1}\in A_{X_{0}}^{0,1}\left( \frak{D}_{1}\left( L_{0}\right) \right)
\end{equation*}
must be $\overline{\partial }$-closed by the integrability conditions in
Lemma \ref{1.6}. Its cohomology class in 
\begin{equation*}
H^{1}\left( \frak{D}_{1}\left( L_{0}\right) \right)
\end{equation*}
is the first-order deformation of the pair $\left( X_{0},L_{0}\right) $
given by $\left( \ref{2.2.1}\right) $ (see \cite{AC}). Its symbol is just
the element of $H^{1}\left( T_{X_{0}}\right) $ giving the Kodaira-Spencer
class for the compatible first-order deformation of the manifold $X_{0}$.

\begin{lemma}
\label{2.4}i) If $X_{0}$ is a Riemann surface $C_{0}$, the space of all
(formal) deformation/trivializations of the pair (curve, line bundle) taken
modulo holomorphic isomorphisms over $\Delta $, is naturally the space of
power series in $t$ with coefficients $\xi _{i}\in A_{C_{0}}^{0,1}\left( 
\frak{D}_{1}\left( L_{0}\right) \right) .$

ii) In general, a (formal) holomorphic section of $L$ is a power series 
\begin{equation*}
s=\sum\nolimits_{i}t^{i}s_{i}
\end{equation*}
with coefficients $s_{i}$ which are $C^{\infty }$-sections of $L_{0}$ such
that 
\begin{equation*}
\sum\nolimits_{i=0}^{\infty }\left( \bar{\partial }s_{i}\right)
t^{i}-\sum\nolimits_{i=0,j=1}^{\infty ,\infty }\xi _{j}\left( s_{i}\right)
t^{i+j}=0.
\end{equation*}

iii) Suppose 
\begin{equation*}
f\in H^{0}\left( L\right) 
\end{equation*}
has divisor $D$ such that 
\begin{equation*}
D_{0}=D\cdot X_{0}
\end{equation*}
is smooth and reduced. Then there is a trivialization 
\begin{equation*}
F_{\sigma }:X\rightarrow X_{0}\times \Delta 
\end{equation*}
such that 
\begin{equation*}
\sigma ^{-1}\left( D_{0}\right) =D,
\end{equation*}
and a unique $\sigma $-compatible trivialization 
\begin{equation*}
F_{\lambda }:L^{\vee }\rightarrow L_{0}^{\vee }\times \Delta 
\end{equation*}
such that 
\begin{equation*}
f=f_{0}\circ \lambda 
\end{equation*}
where 
\begin{equation*}
f_{0}=\left. f\right| _{X_{0}}.
\end{equation*}
We call the trivialization $F_{\lambda }$ \textit{adapted} to the section $f$%
.
\end{lemma}

\begin{proof}
i) By $\left( \ref{2.3.4}\right) $ and Lemma \ref{1.6} all integrability
conditions vanish automatically.

ii) is immediate from Lemma \ref{1.5}.

iii) Let $N$ be a tubular neighborhood of $D_{0}$ in $X$. On $N$ use a
partition-of-unity argument as in \S 5 of \cite{C1} to construct a $%
C^{\infty }$ -retraction 
\begin{equation*}
\upsilon :N\rightarrow N\cap D_{0}
\end{equation*}
such that each fiber is an analytic polydisk. Cover $N$ as above by
coordinate disks $\left\{ W\right\} $. For each $W_{0}=W\cap X_{0}$ which
meets $D_{0}$ construct a \textit{holomorphic} projection 
\begin{equation*}
\upsilon ^{-1}\left( W_{0}\cap D_{0}\right) \rightarrow W_{0}
\end{equation*}
which takes 
\begin{equation*}
\left( W\cap D\right) \rightarrow \left( W_{0}\cap D_{0}\right) .
\end{equation*}
Again as in \S 5 of \cite{C1}, use a C$^{\infty }$-partition-of-unity
argument to ``average'' these local projections to obtain a projection 
\begin{equation*}
\varkappa :N\rightarrow N\cap X_{0}
\end{equation*}
such that 
\begin{equation*}
\upsilon \circ \varkappa =\upsilon
\end{equation*}
and such that, for each $x_{0}\in D_{0}$, 
\begin{equation*}
\left. \varkappa \right| _{\upsilon ^{-1}\left( x_{0}\right) }
\end{equation*}
is holomorphic. $\varkappa $ gives a projection $\sigma $ in some
neighborhood $D$ such that 
\begin{equation*}
D=\sigma ^{-1}\left( D_{0}\right) .
\end{equation*}
Extend by a partition of unity argument to obtain $\sigma :X\rightarrow
X_{0} $ and the corresponding trivialization $F_{\sigma }=\left( \sigma ,\pi
\right) $.

Now let 
\begin{equation*}
L_{0}=\mathcal{O}_{X_{0}}\left( D_{0}\right) .
\end{equation*}
and suppose $D$ is given by local defining functions. Then, on each slice 
\begin{equation*}
\upsilon ^{-1}\left( x_{0}\right) ,
\end{equation*}
$x_{0}\in D_{0}$, the invertible holomorphic functions 
\begin{equation*}
\frac{z_{W}}{z_{W_{0}}\circ \sigma }
\end{equation*}
fit together to give an invertible $C^{\infty }$-function on $W\subset N$ so
that 
\begin{equation*}
h_{W}:=\frac{z_{W}\circ F_{\sigma }^{-1}}{z_{W_{0}}}
\end{equation*}
is an invertible $C^{\infty }$-function on $W_{0}\times \Delta $. If $%
W\nsubseteqq N$ put 
\begin{equation*}
h_{W}=1.
\end{equation*}

So for patching data 
\begin{eqnarray*}
g^{_{WW^{\prime }}}\left( x\right) &=&\sum g_{i}^{_{WW^{\prime }}}\left(
x_{0}\right) t^{i} \\
&=&g^{W_{0}W_{0}^{\prime }}\left( x_{0}\right) \exp \left(
\sum\nolimits_{j>0}a_{j}^{_{WW^{\prime }}}\left( x_{0}\right) t^{j}\right)
\end{eqnarray*}
we have 
\begin{equation*}
\sum\nolimits_{j>0}a_{j}^{_{WW^{\prime }}}\left( x_{0}\right) t^{j}=\log
h_{W^{\prime }}-\log h_{W}
\end{equation*}
The $\sigma $-compatible trivialization $F_{\lambda }$ constructed in $%
\left( \ref{2.3.0}\right) $ is given in this case by 
\begin{equation*}
\left( x,v\right) \mapsto \left( x_{0},\exp \left( \sum\nolimits_{W^{\prime
}}\rho _{W_{0}^{\prime }}\left( x_{0}\right) \left( \log h_{W^{\prime
}}-\log h_{W}\right) \right) \cdot v\right) .
\end{equation*}
So, under this trivialization$,$ $z_{W}$ corresponds to the section of $%
L_{0}^{\vee }\times \Delta $ given over $\left( x_{0},t\right) \in
W_{0}\times \Delta $ by $\left( v,t\right) $ where 
\begin{equation*}
\begin{array}{c}
v=\frac{z_{W_{0}}}{z_{W}\circ F_{\sigma }^{-1}}\exp \left(
\sum\nolimits_{W^{\prime }}\rho _{W_{0}^{\prime }}\left( x_{0}\right) \left(
\log h_{W^{\prime }}\right) \right) \cdot \left( z_{W}\circ F_{\sigma
}^{-1}\right) \\ 
=\exp \left( \sum\nolimits_{W^{\prime }}\rho _{W_{0}^{\prime }}\left(
x_{0}\right) \left( \log h_{W^{\prime }}\right) \right) \cdot z_{W_{0}}.
\end{array}
\end{equation*}
Now replace $\lambda $ with 
\begin{equation*}
\frac{\lambda }{\exp \left( \sum\nolimits_{W^{\prime }}\rho _{W_{0}^{\prime
}}\left( x_{0}\right) \left( \log h_{W^{\prime }}\right) \right) \cdot
\lambda }.
\end{equation*}
\end{proof}

\subsection{Schiffer-type deformations of line bundles\label{2.5}}

We next wish to consider a very special type of line bundle deformation. Our
aim is to be able to apply Lemmas \ref{1.12}-\ref{1.13} to a case in which $%
M_{0}=L_{0}^{\vee }$ is the total space of a line bundle and the holomorphic
functions under consideration are the holomorphic sections of $L_{0}$. Let $%
X_{\beta }/\Delta $ be a Schiffer-type deformation as in \ref{1.8}. That is,
referring to $\left( \ref{1.8.3}\right) $, suppose that $X_{\beta }/\Delta $
is given by Kuranishi data 
\begin{equation*}
\xi _{X_{0}}^{\beta }=\frac{e^{L_{\beta }}-1}{L_{\beta }}\left( \left[ 
\overline{\partial },L_{\beta }\right] \right)
\end{equation*}
for some divisor 
\begin{equation*}
A_{0}\subseteq X_{0}.
\end{equation*}
Let 
\begin{equation*}
A_{\beta }
\end{equation*}
denote the extension of $A_{0}$ to a divisor on $X_{\beta }$.

Let $L/X_{\beta }$ be a deformation of $L_{0}/X_{0}$. By Lemma \ref{1.15}
and $\left( \ref{2.3.0}\right) $ there are compatible trivializations 
\begin{eqnarray*}
F_{\beta } &:&X_{\beta }\rightarrow X_{0}\times \Delta \\
F_{\lambda } &:&L^{\vee }\rightarrow L_{0}^{\vee }\times \Delta .
\end{eqnarray*}
We need that $F_{\lambda }=F_{\tilde{\beta}}$ for some lifting $\tilde{\beta}
$ of $\beta $ to a vector field on $L_{0}^{\vee }\times \Delta $ for which 
\begin{equation*}
\left[ \tilde{\beta},\chi \right] =0.
\end{equation*}

\begin{lemma}
\label{2.6}i) Suppose that $L_{0}$ is trivial over a neighborhood of $A_{0}$
and that the mapping 
\begin{equation*}
H^{0}\left( \frac{\mathcal{O}_{X_{0}}\left( \infty \cdot A_{0}\right) }{%
\mathcal{O}_{X_{0}}}\right) \rightarrow H^{1}\left( \mathcal{O}%
_{X_{0}}\right) 
\end{equation*}
induced by the exact sequence 
\begin{equation*}
0\rightarrow \mathcal{O}_{X_{0}}\rightarrow \mathcal{O}_{X_{0}}\left( \infty 
\cdot A_{0}\right) \rightarrow \frac{\mathcal{O}_{X_{0}}\left( \infty \cdot %
A_{0}\right) }{\mathcal{O}_{X_{0}}}\rightarrow 0
\end{equation*}
is surjective. Then there is a lifting 
\begin{equation*}
\tilde{\beta}
\end{equation*}
of $\beta $ to a vector field on $L_{0}^{\vee }\times \Delta $ which is
meromorphic above 
\begin{equation*}
A_{0}
\end{equation*}
and otherwise $C^{\infty }$ such that $F_{\tilde{\beta}}$ is a
trivialization of $L^{\vee }$.

ii) Referring to i), suppose that 
\begin{equation*}
L=\mathcal{O}_{X_{\beta }}\left( D\right) 
\end{equation*}
and 
\begin{equation*}
\Phi \circ F_{\beta }\left( D\right) =D_{0}\times \Delta 
\end{equation*}
where $D_{0}$ is the zero-scheme associated to a holomorphic section 
\begin{equation*}
f_{0}:L_{0}^{\vee }\rightarrow \Bbb{C}\text{.}
\end{equation*}
Suppose further that $\Phi $ is holomorphic in a neighborhood of $A_{0}%
\times \Delta $. Then there is a lifting $\tilde{\Phi}$ of $\Phi $ so that
the section 
\begin{equation*}
f_{0}\circ \tilde{\Phi}\circ F_{\tilde{\beta}}
\end{equation*}
is a holomorphic section of $L$.
\end{lemma}

\begin{proof}
i) Since $L$ is trivial near $A_{0},$ we can lift $\beta $ to a vector field 
$\tilde{\beta}$ commuting with $\chi $ and meromorphic near $A_{0}$ by a
patching argument as in \ref{2.3}. Any two liftings differ by a vector field 
\begin{equation*}
a\chi =\sum\nolimits_{j>0}a_{j}\chi t^{j}
\end{equation*}
where the $a_{j}$ are fuctions on $X_{0}$ which are meromorphic near $A_{0}$
and $C^{\infty }$ elsewhere. Given that modulo $t^{n}$%
\begin{equation}
L^{\vee }=L_{\tilde{\beta}}^{\vee }  \label{2.6.1}
\end{equation}
we use the surjectivity hypothesis in the statement of the lemma to choose $%
a_{n+1}$ and achieve $\left( \ref{2.6.1}\right) $ modulo $t^{n+1}$.

ii) The deformation $X_{\beta }$ is trivial in a neighborhood of $%
A_{0}\times \Delta $, so we can choose a lifting $\tilde{\Phi}^{\prime }$ of 
$\Phi $ which is holomorphic near $A_{0}\times \Delta $ and extend by a
partition-of-unity argument. Referring to Lemma \ref{2.4}iii), $\tilde{\Phi}%
^{\prime }\circ F_{\tilde{\beta}}=\left( \sigma ^{\prime },\pi \right) $ and
the adapted trivialization $\left( \sigma ,\pi \right) $ are related by 
\begin{equation*}
\sigma =e^{b}\sigma ^{\prime }
\end{equation*}
for some $C^{\infty }$-function $b$ on $X_{0}\times \Delta $. Now set 
\begin{equation*}
\tilde{\Phi}=e^{b}\tilde{\Phi}^{\prime }.
\end{equation*}
\end{proof}

\subsection{Differential operators and basepoint-free systems\label{2.7}}

Suppose now that 
\begin{equation*}
H^{0}\left( L_{0}\right)
\end{equation*}
is basepoint-free. Fix a section

Let $\Bbb{P}_{0}=\Bbb{P}\left( H^{0}\left( L_{0}\right) \right) $ and let 
\begin{eqnarray*}
\nu &:&\Bbb{P}_{0}\times X_{0}\rightarrow \Bbb{P}_{0} \\
\rho &:&\Bbb{P}_{0}\times X_{0}\rightarrow X_{0}
\end{eqnarray*}
be the two projections. Let 
\begin{equation*}
\tilde{L}_{0}\left( 1\right) =\nu ^{*}\mathcal{O}_{\Bbb{P}_{0}}\left(
1\right) \otimes \rho ^{*}L_{0}.
\end{equation*}
Then by the Leray spectral sequence there are natural isomorphisms 
\begin{equation}
\begin{array}{c}
\rho _{*}\tilde{L}_{0}\left( 1\right) =L_{0}\otimes H^{0}\left( L_{0}\right)
^{\vee } \\ 
H^{k}\left( \tilde{L}_{0}\left( 1\right) \right) =H^{k}\left( L_{0}\right)
\otimes H^{0}\left( L_{0}\right) ^{\vee }.
\end{array}
\label{2.7.1}
\end{equation}
There is a tautological section 
\begin{equation}
\tilde{f}_{0}\in H^{0}\left( \tilde{L}_{0}\left( 1\right) \right)
=H^{0}\left( L_{0}\right) \otimes H^{0}\left( L_{0}\right) ^{\vee
}=End\left( H^{0}\left( L_{0}\right) \right)  \label{2.7.2}
\end{equation}
given by the identity map on $H^{0}\left( L_{0}\right) $. Furthermore 
\begin{equation}
\rho _{*}\left( \tilde{f}_{0}\right)  \label{2.7.3}
\end{equation}
is given by the tautological homomorphism 
\begin{equation*}
H^{0}\left( L_{0}\right) \otimes \mathcal{O}_{X_{0}}\rightarrow L_{0}.
\end{equation*}

Also one easily shows by induction using the Euler sequence that 
\begin{equation*}
H^{i}\left( \frak{D}_{n}\left( \mathcal{O}_{\Bbb{P}_{0}}\left( 1\right)
\right) \right) =0
\end{equation*}
for all $i>0$, so also 
\begin{equation*}
R^{i}\rho _{*}\frak{D}_{n}\left( \tilde{L}_{0}\left( 1\right) \right) =0
\end{equation*}
and 
\begin{equation*}
H^{1}\left( \frak{D}_{n}\left( \tilde{L}_{0}\left( 1\right) \right) \right)
=H^{1}\left( \rho _{*}\frak{D}_{n}\left( \tilde{L}_{0}\left( 1\right)
\right) \right) .
\end{equation*}
There is a natural map 
\begin{equation*}
h:\rho _{*}\frak{D}_{n}\left( \tilde{L}_{0}\left( 1\right) \right)
\rightarrow \frak{D}_{n}\left( \rho _{*}\tilde{L}_{0}\left( 1\right) \right)
\end{equation*}
and 
\begin{eqnarray*}
\frak{D}_{n}\left( \rho _{*}\tilde{L}_{0}\left( 1\right) \right) &=&\frak{D}
_{n}\left( L_{0}\otimes H^{0}\left( L_{0}\right) ^{\vee }\right) \\
&=&\frak{D}_{n}\left( L_{0}\right) \otimes End\left( H^{0}\left(
L_{0}\right) \right) .
\end{eqnarray*}

Now via the trace map we have a canonical splitting 
\begin{equation*}
End\left( H^{0}\left( L_{0}\right) \right) =\Bbb{C}\cdot 1\oplus
End^{0}\left( H^{0}\left( L_{0}\right) \right)
\end{equation*}
where $End^{0}$ denotes trace-zero endomorphisms. Notice that 
\begin{equation*}
\frak{D}_{n}^{\prime }:=h\left( \rho _{*}\frak{D}_{n}\left( \tilde{L}%
_{0}\left( 1\right) \right) \right) =\frak{D}_{n}\left( L_{0}\right) \otimes
1\oplus \frak{D}_{n-1}\left( L_{0}\right) \otimes End^{0}\left( H^{0}\left(
L_{0}\right) \right)
\end{equation*}
so we have that 
\begin{equation*}
\frak{D}_{0}^{\prime }=\mathcal{O}_{X_{0}}
\end{equation*}
and we have the exact sequence 
\begin{equation}
0\rightarrow \frak{D}_{n}^{\prime }\rightarrow \frak{D}_{n+1}^{\prime }%
\overset{symbol}{\longrightarrow }\left( S^{n+1}\left( T_{X_{0}}\right)
\otimes 1\right) \oplus \left( S^{n}\left( T_{X_{0}}\right) \otimes
End^{0}\left( H^{0}\left( L_{0}\right) \right) \right) \rightarrow 0
\label{2.7.4}
\end{equation}
is exact.

\subsection{Extendable linear systems on families of curves\label{2.8}}

If $M$ denotes a sufficiently small analytic neighborhood of a general point
in the moduli space of curves of genus $g$, with universal curve $C/M$,
there is a stratification of the locus 
\begin{equation*}
Z_{d}^{r}=\left\{ L:L\,globally\,generated,\,h^{0}\left( L\right)
=r+1\right\} \subseteq Pic^{d}\left( C/M\right)
\end{equation*}
such that all strata are smooth and the projection of each to $M$ is
submersive with diffeomeorphic fibers. Next consider the induced stratification of the pre-image of  $Z_{d}^{r}$ under the Abel-Jacobi map
\begin{equation*}
\alpha :C^{\left( d\right) }/M\rightarrow Pic^{d}\left( C/M\right).
\end{equation*}
By considering the contact locus between this pre-image stratification and the various diagonal loci in $C^{\left( d\right) }/M$, one can construct a refinement of the stratification of 
\begin{equation*}
\alpha ^{-1}\left( Z_{d}^{r}\right) \subseteq C^{\left( d\right) }/M
\end{equation*}
such that all strata are smooth and the projection of each to $M$ is
submersive with diffeomorphic fibers and having the additional property that, beginning with the
initial element $\left( d\right) $ of the partially ordered set $\left\{
\left( d_{1},\ldots ,d_{s}\right) \right\} $ of all partitions of $d$, the
stratification is compatible with each set 
\begin{equation*}
diag_{\left( d_{1},\ldots ,d_{s}\right) }\left( C^{\left( d\right)
}/M_{g}\right) \cap \alpha ^{-1}\left( Z_{d}^{r}\right) .
\end{equation*}

Suppose now that $C_{0}$ is a compact Riemann surface of genus $g$ of
general moduli and that $L_{0}$ is a line bundle of degree $d$ on $C_{0}$
such that the linear system $\Bbb{P}_{0}:=\Bbb{P}\left( H^{0}\left(
L_{0}\right) \right) $ is basepoint-free. Let $C_{\beta }/\Delta $ be a
Schiffer variation supported at a finite set $A_{0}\subseteq C_{0}$. Then,
by genericity of $C_{0}$ and the remarks just above, there is a deformation $%
\Bbb{P}_{\Delta }\subseteq C_{\beta }^{\left( d\right) }$ over $\Delta $ of $%
\Bbb{P}_{0}\subseteq C_{0}^{\left( d\right) }$ for which there exists a
trivialization 
\begin{equation}
T:\Bbb{P}_{\Delta }\rightarrow \Bbb{P}_{0}\times \Delta  \label{3.5.1}
\end{equation}
compatible with each partition locus of $d$, that is, for each partition $%
\left( d_{1},\ldots ,d_{s}\right) $ of $d$, 
\begin{equation}
T\left( diag_{\left( d_{1},\ldots ,d_{s}\right) }\left( C_{\beta }^{\left(
d\right) }\right) \times _{C_{\beta }^{\left( d\right) }}\Bbb{P}_{\Delta
}\right) =\left( diag_{\left( d_{1},\ldots ,d_{s}\right) }\left(
C_{0}^{\left( d\right) }\right) \times _{C_{0}^{\left( d\right) }}\Bbb{P}%
_{0}\right) \times \Delta .  \label{3.5.2}
\end{equation}
Notice that $T$ is a $C^{\infty }$-map, and is not in general analytic.
However $T$ can be chosen so that, for each $p\in \Bbb{P}_{0}$, $%
T^{-1}\left( \left\{ p\right\} \times \Delta \right) $ is a proper analytic
subvariety of $\Bbb{P}_{\Delta }$.

Now the tautological section $\tilde{f}_{0}$ of $\tilde{L}_{0}\left(
1\right) =\mathcal{O}_{\Bbb{P}_{0}}\boxtimes L_{0}$ defined in $\left( \ref
{2.7.2}\right) $ has divisor 
\begin{equation*}
D_{0}\subseteq \Bbb{P}_{0}\times C_{0}.
\end{equation*}
Let 
\begin{equation*}
D\subseteq \Bbb{P}_{\Delta }\times _{\Delta }C_{\beta }
\end{equation*}
denote the divisor of the tautological section $\tilde{f}$ of 
\begin{equation*}
\tilde{L}\left( 1\right) :=\mathcal{O}_{\Bbb{P}_{\Delta }}\left( 1\right)
\boxtimes _{\Delta }L.
\end{equation*}
Then, by $\left( \ref{3.5.2}\right) ,$ the ``product'' trivialization 
\begin{equation*}
\left( T,F_{\beta }\right) :\Bbb{P}_{\Delta }\times _{\Delta }C_{\beta
}\rightarrow \Bbb{P}_{0}\times C_{0}\times \Delta 
\end{equation*}
is compatible with the trivialization $T$ in $\left( \ref{3.5.1}\right) $,
that is, for each $p\in \Bbb{P}_{0}$, 
\begin{equation*}
\left( T,F_{\beta }\right) ^{-1}\left( \left\{ p\right\} \times C_{0}\times 
\Delta \right) =T^{-1}\left( \left\{ p\right\} \times \Delta \right) \times %
_{\Bbb{P}_{\Delta }}\left( \Bbb{P}_{\Delta }\times _{\Delta }C_{\beta
}\right) .
\end{equation*}
That is, we have the commutative diagram
\begin{equation*}
\begin{array}{ccc}
\Bbb{P}_{\Delta }\times _{\Delta }C_{\beta } & \overset{\left( T,F_{\beta
}\right) }{\longrightarrow } & \Bbb{P}_{0}\times C_{0}\times \Delta  \\ 
\downarrow  &  & \downarrow  \\ 
\Bbb{P}_{\Delta } & \overset{T}{\longrightarrow } & \Bbb{P}_{0}\times \Delta 
\end{array}
\end{equation*}

Furthermore, by $\left( \ref{3.5.2}\right) $, we can adjust $\left(
T,F_{\beta }\right) $ ``in the $C_{0}$-direction'' to obtain a trivialization
\begin{equation*}
\begin{array}{ccc}
\Bbb{P}_{\Delta }\times _{\Delta }C_{\beta } & \overset{F}{\longrightarrow }
& \Bbb{P}_{0}\times C_{0}\times \Delta  \\ 
\downarrow  &  & \downarrow  \\ 
\Bbb{P}_{\Delta } & \overset{T}{\longrightarrow } & \Bbb{P}_{0}\times \Delta 
\end{array}
\end{equation*}
which maintains the property 
\begin{equation}
F^{-1}\left( \left\{ p\right\} \times C_{0}\times \Delta \right)
=T^{-1}\left( \left\{ p\right\} \times \Delta \right) \times _{\Bbb{P}_{%
\Delta }}\left( \Bbb{P}_{\Delta }\times _{\Delta }C_{\beta }\right) .
\label{3.5.3}
\end{equation}
and achieves in addition that 
\begin{equation}
F^{-1}\left( D_{0}\times \Delta \right) =D.  \label{3.5.4}
\end{equation}
Finally, we can choose the adjustments to be holomorphic in the $C_{0}$%
-direction in a small neighborhood of $\Bbb{P}_{\Delta }\times A_{0}\times 
\Delta $.

Thus referring to Lemma \ref{1.12} there is a $C^{\infty }$-vector field 
\begin{equation*}
\gamma =\sum\nolimits_{n>0}\gamma _{n}t^{n}
\end{equation*}
on $\Bbb{P}_{0}\times C_{0}\times \Delta $ of type $\left( 1,0\right) $ such
that

1) each $\gamma _{n}$ annihilates functions pulled back from $\Bbb{P}_{0}$,
that is, it is an $\mathcal{O}_{\Bbb{P}_{0}}$-linear operator,

2) for each $n$ and each $p\in \Bbb{P}_{0}$, 
\begin{equation*}
\left. \gamma _{n}\right| _{\left\{ p\right\} \times C_{0}}
\end{equation*}
is meromorphic on a neighborhood of $\left\{ p\right\} \times A_{0}$,

3) given a function 
\begin{equation*}
g=\sum\nolimits_{k=0}^{\infty }g_{k}t^{k}:\Bbb{P}_{0}\times C_{0}\times
\Delta \rightarrow \Bbb{C}
\end{equation*}
with each $g_{k}$ a $C^{\infty }$-function on (an open set in) $\Bbb{P}%
_{0}\times C_{0}$ and any point $p\in \Bbb{P}_{0}$, 
\begin{equation*}
\left. g\circ F\right| _{T^{-1}\left( \left\{ p\right\} \times \Delta
\right) \times _{\Bbb{P}_{\Delta }}\left( \Bbb{P}_{\Delta }\times _{\Delta
}C_{\beta }\right) }
\end{equation*}
is holomorphic if and only if 
\begin{equation*}
\left. \left[ \overline{\partial _{0}},e^{L_{-\gamma }}\right] \left(
g\right) \right| _{\left\{ p\right\} \times C_{0}\times \Delta }=0.
\end{equation*}

Again, following Lemma \ref{2.6}, there is a trivalization 
\begin{equation*}
\begin{array}{ccc}
\tilde{L}\left( 1\right) ^{\vee } & \overset{\tilde{F}}{\longrightarrow } & 
\tilde{L}_{0}\left( 1\right) ^{\vee }\times \Delta  \\ 
\downarrow  &  & \downarrow  \\ 
\Bbb{P}_{\Delta }\times _{\Delta }C_{\beta } & \overset{F}{\longrightarrow }
& \Bbb{P}_{0}\times C_{0}\times \Delta 
\end{array}
\end{equation*}
of $\tilde{L}\left( 1\right) $ and a lifting $\tilde{\gamma}$ of $\gamma $
such that, for the tautological sections $\tilde{f}_{0}$ and $\tilde{f}$
defined earlier in this section, 
\begin{equation*}
\tilde{f}=\tilde{F}\circ \tilde{f}_{0}.
\end{equation*}
Thus, for each $p\in \Bbb{P}_{0}$, 
\begin{equation}
\left. \left[ \overline{\partial _{0}},e^{L_{-\tilde{\gamma}}}\right] \left( 
\tilde{f}_{0}\right) \right| _{\left\{ p\right\} \times C_{0}\times \Delta %
}=0.  \label{3.5.5}
\end{equation}
Let 
\begin{equation*}
\frak{D}_{n}^{\Bbb{P}_{0}}\left( \tilde{L}_{0}\left( 1\right) \right)
\subseteq \frak{D}_{n}\left( \tilde{L}_{0}\left( 1\right) \right) 
\end{equation*}
denotes the subsheaf of $\mathcal{O}_{\Bbb{P}_{0}}$-linear operators. Then 
\begin{equation*}
\left[ \overline{\partial _{0}},e^{L_{-\tilde{\gamma}}}\right] 
\end{equation*}
is a $\overline{\partial _{0}}$-closed element of 
\begin{equation*}
\sum\nolimits_{n>0}^{\infty }H^{1}\left( \frak{D}_{n}^{\Bbb{P}_{0}}\left( 
\tilde{L}_{0}\left( 1\right) \right) \right) t^{n}.
\end{equation*}

Now, referring to $\left( \ref{2.7.4}\right) $, we need to analyze 
\begin{equation*}
\rho _{*}\left[ \overline{\partial _{0}},e^{L_{-\tilde{\gamma}}}\right] \in
\sum\nolimits_{n>0}H^{1}\left( \frak{D}_{n}^{\prime }\right)
t^{n}=\sum\nolimits_{n>0}H^{1}\left( \frak{D}_{n}\left( L_{0}\right) \right)
\otimes End\left( H^{0}\left( L_{0}\right) \right) t^{n}.
\end{equation*}
In fact, by construction, this element lies in the image of 
\begin{equation*}
\sum\nolimits_{n>0}^{\infty }H^{1}\left( \rho _{*}\frak{D}_{n}^{\Bbb{P}%
_{0}}\left( \tilde{L}\left( 1\right) \right) \right)
t^{n}=\sum\nolimits_{n>0}H^{1}\left( \frak{D}_{n}\left( L_{0}\right) \right)
\otimes \Bbb{C}\cdot \left( id\right) \cdot t^{n}\subseteq
\sum\nolimits_{n>0}H^{1}\left( \frak{D}_{n}\left( L_{0}\right) \right)
\otimes End\left( H^{0}\left( L_{0}\right) \right) t^{n}.
\end{equation*}
Now 
\begin{equation*}
H^{1}\left( \tilde{L}_{0}\left( 1\right) \right) =Hom\left( H^{0}\left(
L_{0}\right) ,H^{1}\left( L_{0}\right) \right) .
\end{equation*}
But by $\left( \ref{3.5.5}\right) $, the image of 
\begin{equation*}
\left. \left\{ \left[ \overline{\partial _{0}},e^{L_{-\tilde{\gamma}%
}}\right] \left( \tilde{f}_{0}\right) \right\} \right| _{\left\{ p\right\} 
\times C_{0}\times \Delta }\in \sum\nolimits_{n>0}^{\infty }H^{1}\left(
L_{0}\right) \cdot t^{n}.
\end{equation*}
is zero for each $p\in \Bbb{P}_{0}$. Thus 
\begin{equation}
\rho _{*}\left[ \overline{\partial _{0}},e^{L_{-\tilde{\gamma}}}\right]
\left( \rho _{*}\tilde{f}_{0}\right) =0\in \sum\nolimits_{n>0}Hom\left(
H^{0}\left( L_{0}\right) ,H^{1}\left( L_{0}\right) \right) t^{n}.
\label{2.8.6}
\end{equation}

\begin{theorem}
\label{2.9}Suppose $X_{0}$ is a curve of genus $g$ of general moduli.
Suppose further that, by varying of $\beta $ in \ref{2.5}, the coefficients
to $t^{n+1}$ in all expressions 
\begin{equation*}
\left[ \overline{\partial },e^{-L_{\beta }}\right] 
\end{equation*}
generate $H^{1}\left( S^{n+1}\left( T_{X_{0}}\right) \right) $ for each $n%
\geq 0$. (For example we allow the divisor $A_{0}\subseteq X_{0}$ to move.)
Then the maps 
\begin{equation*}
\mu ^{n+1}:H^{1}\left( S^{n+1}T_{X_{0}}\right) \rightarrow \frac{Hom\left(
H^{0}\left( L_{0}\right) ,H^{1}\left( L_{0}\right) \right) }{image\ \tilde{%
\mu}^{n}}
\end{equation*}
are zero for all $n\geq 0$.
\end{theorem}

\begin{proof}
Let 
\begin{equation*}
\rho _{*}\left[ \overline{\partial _{0}},e^{L_{-\tilde{\gamma}}}\right]
_{n+1}.
\end{equation*}
denote the coefficient of $t^{n+1}$ in $\rho _{*}\left[ \overline{\partial %
_{0}},e^{L_{-\tilde{\gamma}}}\right] $. Referring to $\left( \ref{2.7.4}%
\right) $ and the fact the the operators take values in the sheaf $\frak{D}%
_{n}^{\Bbb{P}_{0}}\left( \tilde{L}_{0}\left( 1\right) \right) $, we have
that 
\begin{equation}
symbol\left( \left( \rho _{*}\left[ \overline{\partial _{0}},e^{L_{-\tilde{%
\gamma}}}\right] \right) _{n+1}\right) =\left( \overline{\partial }\beta
_{1}^{n+1}\otimes 1\right) \oplus 0\in S^{n+1}\left( T_{X_{0}}\right) \oplus
\left( S^{n}\left( T_{X_{0}}\right) \otimes End^{0}\left( H^{0}\left(
L_{0}\right) \right) \right)   \label{2.10.1}
\end{equation}
where 
\begin{equation*}
\beta =\sum\nolimits_{j>0}\beta _{j}t^{j}.
\end{equation*}
By $\left( \ref{2.10.1}\right) $ and the hypothesis that the elements $%
\overline{\partial }\beta _{1}^{n+1}$ generate $H^{1}\left(
S^{n+1}T_{X_{0}}\right) $, we have that, by varying $\beta $, the elements 
\begin{equation*}
symbol\left( \rho _{*}\left[ \overline{\partial _{0}},e^{L_{-\tilde{\gamma}%
}}\right] _{n+1}\right) 
\end{equation*}
generate 
\begin{equation*}
S^{n+1}\left( T_{X_{0}}\right) 
\end{equation*}
for each $n\geq 0$.

Thus, by $\left( \ref{2.7.4}\right) $ and $\left( \ref{2.8.6}\right) $, the
map $\tilde{\nu}^{n+1}$ given by 
\begin{eqnarray*}
H^{1}\left( \frak{D}_{n+1}\left( L_{0}\right) \right) &\rightarrow &\frac{%
Hom\left( H^{0}\left( L_{0}\right) ,H^{1}\left( L_{0}\right) \right) }{%
image\left( \tilde{\nu}^{n}\right) } \\
D &\mapsto &D\left( \tilde{f}_{0}\right)
\end{eqnarray*}
is zero for all $n\geq 0$.
\end{proof}

\section{Brill-Noether theory\label{3}}

In this last section we give a simple application of Theorem \ref{2.9} to
Brill-Noether theory. From now on we assume that $X_{0}$ is a \textit{generic%
} compact Riemann surface $C_{0}$. We choose 
\begin{equation*}
A_{0}=\left\{ x_{0}\right\}
\end{equation*}
in \S \ref{1}-\ref{2} where $x_{0}$ is a general point of $C_{0}$, and let 
\begin{equation*}
C_{\beta }/\Delta
\end{equation*}
denote the family of Schiffer-type deformation associated to some vector
field 
\begin{equation*}
\beta =\sum\nolimits_{j>0}\beta _{j}t^{j}
\end{equation*}
where each $\beta _{j}$ is meromorphic with poles in some neighborhood $%
U_{0} $ of $x_{0}$. Since $C_{0}$ is generic, there exists a line-bundle
deformation 
\begin{equation*}
L/C_{\beta }
\end{equation*}
such that 
\begin{equation*}
H^{0}\left( L\right) \rightarrow H^{0}\left( L_{0}\right)
\end{equation*}
is surjective. We wish to apply Theorem \ref{2.9} to conclude that the maps $%
\mu ^{n+1}$ are all zero for $n\geq 0$.

\begin{lemma}
\label{3.2}Let $\beta _{1}$ range over all vector fields such that the
Kodaira-Spencer class 
\begin{equation*}
\overline{\partial }\beta _{1}
\end{equation*}
generates the kernel of the map 
\begin{equation*}
H^{1}\left( T_{C_{0}}\right) \rightarrow H^{1}\left( T_{C_{0}}\left(
x_{0}\right) \right) .
\end{equation*}
Then the elements 
\begin{equation*}
\overline{\partial }\left( \beta _{1}^{k+1}\right) 
\end{equation*}
generate the kernel of the map 
\begin{equation*}
H^{1}\left( S^{k+1}T_{C_{0}}\right) \rightarrow H^{1}\left(
S^{k+1}T_{C_{0}}\left( \left( k+1\right) x_{0}\right) \right) .
\end{equation*}
\smallskip 
\end{lemma}

\begin{proof}
Let $z$ be a local analytic coordinate for $C_{0}$ centered on $x_{0}$. We
trivialize our Schiffer-type variation of $C_{0}$ so that 
\begin{equation*}
\beta _{1}=\frac{\rho }{z}\frac{\partial }{\partial z}
\end{equation*}
where $\rho $ is a $C^{\infty }$-function on $C_{0}$ such that

i) $\rho $ is supported on an arbitrarily small neighborhood of $x_{0}$,

ii) in a smaller neighborhood $U_{0}$ of $x_{0}$, 
\begin{equation*}
\rho =\frac{a_{-1}}{z}+a_{0}+\ldots +a_{k}z^{k}.
\end{equation*}
So 
\begin{equation*}
\bar{\partial }\left( \frac{\rho }{z}\right) ^{k+1}\left( \frac{\partial }{%
\partial z}\right) ^{k+1}
\end{equation*}
represents the symbol of $\left[ \overline{\partial },L_{\beta
_{1}}^{k+1}\right] $. By varying the choice of the $a_{i}$ in the definition
of $\rho $ we can therefore obtain symbols which generate the image of 
\begin{equation*}
\frac{S^{k+1}T_{C_{0}}\left( \left( k+1\right) x_{0}\right) }{%
S^{k+1}T_{C_{0}}}
\end{equation*}
in $H^{1}\left( S^{k+1}T_{C_{0}}\right) .$
\end{proof}

Now if $x_{0}$ varies over a dense subset of $C_{0}$, the elements of kernel
of 
\begin{equation*}
H^{1}\left( S^{k+1}T_{C_{0}}\right) \rightarrow H^{1}\left(
S^{k+1}T_{C_{0}}\left( x_{0}\right) \right)
\end{equation*}
generate $H^{1}\left( S^{k+1}T_{C_{0}}\right) $. So we conclude by Theorem 
\ref{2.9}:

\begin{theorem}
\label{3.3}If $C_{0}$ is a curve of general moduli and $H^{0}\left(
L_{0}\right) $ is basepoint-free, the mapping 
\begin{equation*}
\mu ^{k+1}:H^{1}\left( S^{k+1}T_{C_{0}}\right) \rightarrow \frac{Hom\left(
H^{0}\left( L_{0}\right) ,H^{1}\left( L_{0}\right) \right) }{%
\sum\nolimits_{k^{\prime }\leq k}image\ \tilde{\mu}^{k^{\prime }}}
\end{equation*}
given in \ref{2.1} must be the zero map for $k\geq 0$.
\end{theorem}

To see that Petri's conjecture follows from Theorem \ref{3.3}, we reason as
in \S 9 of [ACGH]. Namely we consider the dual mappings 
\begin{equation*}
\mu _{k}:\ker \mu _{k-1}\rightarrow H^{0}\left( \omega _{C_{0}}^{k+1}\right)
\end{equation*}
(inductively defined beginning with the zero map 
\begin{equation*}
\mu _{-1}:H^{0}\left( L_{0}\right) \otimes H^{0}\left( \omega
_{C_{0}}\otimes L_{0}^{\vee }\right) \rightarrow \left\{ 0\right\} ).
\end{equation*}
Petri's conjecture asserts that, for our $C_{0}$ of general moduli, the
mapping 
\begin{equation*}
\mu _{0}:H^{0}\left( L_{0}\right) \otimes H^{0}\left( \omega _{C_{0}}\otimes
L_{0}^{\vee }\right) \rightarrow H^{0}\left( \omega _{C_{0}}\right) ,
\end{equation*}
which is of course simply the multiplication map, is injective. To see that
this follows from Theorem \ref{3.3}, let $\left\{ s_{i}\right\} $ denote a
basis for $H^{0}\left( L_{0}\right) $. Suppose now that 
\begin{equation*}
\mu _{0}\left( \sum s_{i}\otimes t_{i}\right) =0.
\end{equation*}
Then the element 
\begin{equation*}
\sum (ds_{i})t_{i}\in H^{0}\left( \omega _{C_{0}}^{2}\right)
\end{equation*}
is well-defined, giving the mapping $\mu _{1}$, etc. Since, by Theorem \ref
{3.3}, successive maps $\mu _{k}$ are the zero map we have, for any local
trivialization of $L_{0}$ and local coordinate $z$ near a general point $%
x_{0}$ on $C_{0}$, the local system of (pointwise) equations 
\begin{equation*}
\sum_{i}t_{i}\left( x_{0}\right) \frac{d^{k}s_{i}}{dz^{k}}\left(
x_{0}\right) =0
\end{equation*}
for all $k$, which is clearly impossible unless all the $t_{i}\left(
x_{0}\right) $ are zero.

\end{document}